\documentclass[11pt, letterpaper]{amsart}
\usepackage{amsmath}
\usepackage{amssymb}
\usepackage[left=1in,right=1in,bottom=0.8in,top=0.8in]{geometry}
\usepackage{amsfonts}
\usepackage{graphicx}
\usepackage[font=small,labelfont=bf]{caption}
\usepackage{amsthm}
\usepackage[T1]{fontenc}
\usepackage{mathrsfs}
\usepackage{mathtools}
\usepackage{stmaryrd}
\SetSymbolFont{stmry}{bold}{U}{stmry}{m}{n}
\usepackage{enumerate}
\usepackage[dvipsnames]{xcolor}
\usepackage{tikz}
\usepackage{tikz-cd}
\usetikzlibrary{decorations.pathmorphing}
\usepackage{rotating}
\usepackage{url}
\usepackage{hyperref}
\usepackage[
backend=biber,
style=alphabetic,
maxbibnames=99,
]{biblatex}

\makeatletter
\renewcommand{\tocsection}[3]{%
  \indentlabel{\@ifnotempty{#2}{\bfseries\ignorespaces#1 #2.\quad}}\bfseries#3}
  \renewcommand{\tocsubsection}[3]{%
  \indentlabel{\hspace{2.5em}\@ifnotempty{#2}{\ignorespaces#1 #2.\quad}}#3}
\makeatother

\newtheorem{theorem}{Theorem}[section]
\newtheorem{corollary}[theorem]{Corollary}
\newtheorem{proposition}[theorem]{Proposition}
\newtheorem{lemma}[theorem]{Lemma}

\newtheorem*{proposition*}{Proposition}

\theoremstyle{definition}
\newtheorem{definition}[theorem]{Definition}

\theoremstyle{remark}
\newtheorem{remark}[theorem]{Remark}

\newtheorem{innercustomcon}{Conjecture}
\newenvironment{customcon}[1]
  {\renewcommand\theinnercustomcon{#1}\innercustomcon}
  {\endinnercustomcon}

\definecolor{refblue}{RGB}{30,90,180}
\definecolor{signal}{RGB}{203,23,206}
\hypersetup{
    colorlinks=true,
    linkcolor=refblue,
    citecolor=signal,
    urlcolor=signal,
    pdftitle={A conditional proof of the Casas-Alvero conjecture},
    pdfauthor={Soham Ghosh},
    pdfsubject={A conditional proof of the Casas-Alvero conjecture},
    pdfkeywords={Casas-Alvero, Koszul homology, regular sequences, uniform local estimates, elementary symmetric polynomials}
}
\numberwithin{equation}{section}
\makeatletter
\@namedef{subjclassname@2020}{%
  \textup{2020} Mathematics Subject Classification}
\makeatother

\DeclareMathOperator{\Res}{\ensuremath{Res}}

\theoremstyle{definition}
\newtheorem{innercustomhyp}{Hypothesis}
\newenvironment{customhyp}[1]
  {\renewcommand\theinnercustomhyp{#1}\innercustomhyp}
  {\endinnercustomhyp}

\author[Soham Ghosh]{Soham Ghosh}
\address{Department of Mathematics, University of Washington, Seattle, WA 98195, USA}
\email{soham13@uw.edu}
\title{A conditional proof of the Casas-Alvero conjecture}
\date{\today}

\begin{document}
\mathtoolsset{showonlyrefs}

\begin{abstract}
The Casas-Alvero conjecture states that if $f(X)$ is a monic univariate polynomial of degree $d$ over a characteristic $0$ field $\mathbb{K}$ such that $\gcd(f, f_{i})$ is non-trivial for each $i=1, \dots, d-1$, then $f(X)=(X-\alpha)^d$ for some $\alpha\in \mathbb{K}$. In this paper, we prove the Casas-Alvero conjecture for polynomials of any degree $d\geq 3$ over any characteristic $0$ field, under a uniform local {\L}ojasiewicz type  hypothesis for a family of polynomials obtained from the elementary symmetric functions. Furthermore, we prove this estimate on compact subsets of a common Zariski open set of coefficients parametrizing the family.
\end{abstract}
\maketitle

\section{Introduction}\label{sec1}
\let\oldthetheorem\thetheorem

\subsection{Aim and main results of the paper}\label{subsec1.1}
Let $\mathbb{K}$ be a field and $f(X)\in \mathbb{K}[X]$ be a monic polynomial of degree $d>1$. Let $f^{(i)}(X):=d^if(X)/dX^i$ be the $i^{th}$ formal derivative of $f(X)$ with respect to $X$ and let $f_i(X)$ be the $i^{th}$ Hasse--Schmidt derivative of $f(X)$. Over fields of characteristic $0$, the two derivatives are related via $f_i(X)=f^{(i)}(X)/i!$. This paper is concerned with the following question posed by E. Casas-Alvero in connection with his work \cite{CA} on higher-order polar germs:

\begin{customcon}{CA}[Casas-Alvero, 2001]\label{con1}
 Let $f(X)$ be a monic univariate polynomial of degree $d$ over a field $\mathbb{K}$. Then $\gcd(f, f_{i})$ is non-trivial for each $i=1, \dots, d-1$ if and only if $f(X)=(X-\alpha)^d$ for some $\alpha\in \mathbb{K}$.
\end{customcon}

Conjecture~\ref{con1} was originally posed over characteristic $0$ fields using ordinary formal derivatives. Using Hasse-Schmidt derivatives, the conjecture has been studied in the literature in the above generalized form over any field irrespective of the characteristic. We review the existing literature on this problem in Section~\ref{subsec1.2} below. It is well-known that the conjecture is false in general in any positive characteristic, leaving only the characteristic zero case open. A bad prime for that degree is a prime $p$ for which the assertion fails over $\overline{\mathbb F}_p$. Our aim is to prove Conjecture~\ref{con1} in characteristic $0$ under the following hypothesis.

Fix $n\ge3$, a word $\mathbf j=(j_1,\ldots,j_{n-1})$ with $j_i\in[n]:=\{1,2,\dots, n\}$, and an index $u\in U_\mathbf{j}:=[n]\setminus\{j_1,\dots, j_{n-1}\}$. Let $q=n-2$, $z_{j_1}=z_u=0$ and let $\mathbf{y}$ be the tuple of the remaining $q$ coordinates $z_i$ for $1\leq i\leq n$, $i\neq j_1, u$. For $\mathbf c=(c_2,\ldots,c_{n-1})\in\mathbb C^q$, consider the polynomial map $Q_{u,\mathbf c}=(Q_{i,u,\mathbf c})_{i=2}^{n-1}:\mathbb C^q\to\mathbb C^q$ given by
\[
 \begin{aligned}
 Q_{i,u,\mathbf c}(\mathbf{y})=(1-c_i z_{j_i})e_{n-i}\bigl((z_k-z_{j_i})_{k=1}^n\bigr)
        +c_i e_{n+1-i}\bigl((z_k-z_{j_i})_{k=1}^n\bigr).
 \end{aligned}
\]
Here $e_h$ denotes the $h^{th}$ elementary symmetric polynomial. Write $\mathbf1=(1,\ldots,1)\in\mathbb C^q$. Also consider the following finite set of rational points in $\mathbb C^q$:
\[\mathscr G_n=\prod_{i=2}^{n-1}
         \left\{\frac{n+1-i}{r}:1\le r\le i-1\right\}
         \subset\mathbb Q_{>0}^{\,q}.\]

\begin{customhyp}{L}\label{Hypothesis:L}
For every $n\ge4$, every word $\mathbf j$ and every $u\in U_\mathbf{j}$ as above, if $\mathbf{0}$ is an isolated zero of $Q_{u,\mathbf1}$, then there exist $A,B,\varepsilon>0$ such that
\begin{equation}\label{eq:finite-grid-estimate}
 \|Q_{u,\mathbf c}(\mathbf{y})\|\ge A\|\mathbf{y}\|^B
 \qquad\bigl(\|\mathbf{y}\|<\varepsilon,\ 
              \operatorname{dist}(\mathbf c,\mathscr G_n)<\varepsilon\bigr).
\end{equation}
\end{customhyp}
The norms and distances in the hypothesis above are Euclidean. The constants $A, B, \varepsilon$ may depend on $(n,\mathbf j,u)$, but the same constants must work for all coefficients in the indicated neighborhood of $\mathscr G_n$ in $\mathbb{C}^q$. Hypothesis~\ref{Hypothesis:L} can be interpreted as a uniform form of the local {\L}ojasiewicz inequality \cite{Lojasiewicz1959}. For   holomorphic map $F:\mathbb{C}^q\rightarrow\mathbb{C}^q$ with an isolated zero at $\mathbf{0}$,
this inequality gives $\|F(\mathbf{y})\|\ge A\|\mathbf{y}\|^B$
near $\mathbf{0}$ for some $A,B>0$. The map $Q_{u,\mathbf c}:\mathbb{C}^q\rightarrow \mathbb{C}^q$ in Hypothesis~\ref{Hypothesis:L} varies with $\mathbf c\in \mathbb{C}^q$ and for each fixed pair $(\mathbf j,u)$, the same
constants $A,B$ and the same radius must work for all $\mathbf c$
in a neighborhood of $\mathscr G_n$. In the terminology of
\cite[\S6, Definition~6.1 and Remark~6.2]{LejeuneJalabertTeissier2008},
this asks that the {\L}ojasiewicz exponent of $(y_1,\ldots,y_q)$
relative to $(Q_{2,u,\mathbf c},\ldots,Q_{n-1,u,\mathbf c})$ be
finite along $\{\mathbf{0}\}\times\mathscr G_n$ in the space of
pairs $(\mathbf{y},\mathbf c)$. The assertion of
Hypothesis~\ref{Hypothesis:L} is that isolatedness of the origin
for the single map $Q_{u,\mathbf1}$ suffices for this uniform estimate.

The main result of this paper is the following.

\begingroup
\renewcommand*{\thetheorem}{A}
\renewcommand*{\theHtheorem}{main.A}
\begin{theorem}\label{MainConditional}
Assuming Hypothesis~\ref{Hypothesis:L}, Conjecture~\ref{con1} holds over every field of characteristic zero.
\end{theorem}
\endgroup

Our second result concerns the estimate in \eqref{eq:finite-grid-estimate} at general coefficients $\mathbf{c}\in\mathbb{C}^q$.

\begingroup
\renewcommand*{\thetheorem}{B}
\renewcommand*{\theHtheorem}{main.B}
\begin{theorem}\label{MainGeneric}
For every $n\ge3$ there exists a nonzero $\Delta_n\in\mathbb Q[c_2,\ldots,c_{n-1}]$ such that the following hold for every word $\mathbf j$ and every $u\in U_\mathbf{j}$.
\begin{enumerate}[(i)]
\item For each $\mathbf c\in\Omega_n:=\{\mathbf c\in\mathbb C^{n-2}:\Delta_n(\mathbf c)\ne0\}$, 
\[\text{the algebra }
\mathbb C[\mathbf{y}]/(Q_{2,u,\mathbf c},\ldots,Q_{n-1,u,\mathbf c}) \text{ is finite-dimensional.}\]
\item For every compact $K\subset\Omega_n$ there are $A_K,R_K>0$, independent of $\mathbf j$ and $u$, such that
\[
 \|Q_{u,\mathbf c}(\mathbf{y})\|\ge A_K\|\mathbf{y}\|^{(n-1)!}
 \qquad(\mathbf c\in K,\ \|\mathbf{y}\|<R_K).
\]
\end{enumerate}
Moreover, $\Omega_n\cap\mathbb Q_{>0}^{\,n-2}\ne\varnothing$, and
$Q_{2,u,\mathbf c},\ldots,Q_{n-1,u,\mathbf c}$ form a regular sequence in the polynomial ring $\mathbb C[\mathbf{y},\mathbf c]$ before specialization of the coefficients.
\end{theorem}
\endgroup
\setcounter{theorem}{0}

Theorem~\ref{MainGeneric} follows from Proposition~\ref{prop:gw-universal-regularity} and Theorem~\ref{thm:gw-generic-estimate}. In Theorem~\ref{thm:grid-optimal-open} we show existence of $\Delta_n\in\mathbb Q[c_2,\ldots,c_{n-1}]$, whose nonvanishing locus contains all the points of $\mathscr G_n$ that can belong to an open set satisfying both \textup{(i)} and \textup{(ii)}.

\subsubsection{Idea of the proof}
We briefly sketch the proof of Theorem~\ref{MainConditional}. By \cite[Proposition~5.2]{SG}, the conjecture in degree $n$ over any algebraically closed field $\mathbb{K}$ is equivalent to the regularity of a sequence of $n-1$ homogeneous polynomials $f_1,\ldots, f_{n-1}$ in $R_{n-1}:=\mathbb{K}[x_1,\dots, x_{n-1}]$, for every choice of $\mathbf j=(j_1,\dots, j_{n-1})\in[n]^{n-1}$. Let $C(n,w)$ denote this assertion restricted to words $\mathbf{j}$ using at most $w$ distinct indices. We prove the implication
\[
 C(n+1,w)\ \text{and}\ C(n,w-1)\quad\rightarrow\quad C(n,w)
 \qquad(3\le w\le n-1).
\]

Fix a word $\mathbf{j}\in [n]^{n-1}$ using exactly $w$ indices. Adjoining the variable $x_n$ to $R_{n-1}$ relates the generators $f_i$ for $C(n,w)$ to the generators $F_i$ in the next degree $C(n+1,w)$. Filtering the Koszul complex on $F_1,\ldots,F_{n-1}$ by degree in the new variable $x_n$ yields a connecting map from the first Koszul homology of $f_1,\ldots,f_{n-1}$. Under the two induction premises $C(n+1, w)$ and $C(n,w-1)$, the connecting map identifies this homology with the kernel of the map on zeroth homology induced by the first filtration inclusion. On the other hand, the assumption $C(n,w-1)$ implies that $f_2,\ldots,f_{n-1}$ form a regular sequence, so the same homology is the annihilator of $f_1$ in the quotient by $(f_2,\ldots,f_{n-1})$. We show that the factors of $f_1$ corresponding to selected indices in $\mathbf{j}$ are already non-zero divisors and thus, by symmetry it remains to treat $z_u$ for one unused index $u\in U_\mathbf{j}:=[n]\setminus\{j_1,\dots, j_{n-1}\}$.

The induction hypothesis $C(n+1,w)$ yields the isolated zero of $Q_{u,\mathbf1}$ required in Hypothesis~\ref{Hypothesis:L}. A nonzero common zero of the forms $f_i|_{z_u=0}$, $2\le i\le n-1$, would produce, by a fractional linear change of coordinates, nonzero zeros of $Q_{u,\mathbf c}$ tending to the origin with $\mathbf c$ tending to $\mathscr G_n$. This contradicts \eqref{eq:finite-grid-estimate}. The restricted forms are therefore a homogeneous system of parameters, and the Koszul exact sequence for multiplication by $z_u$ gives the required vanishing of the annihilator.

Finally, the result of Draisma and Knopper \cite[Proposition~6]{DJ} supplies the statement $C(d,2)$ for all $d\ge3$. We first induct on $w$, assuming $C(d,w-1)$ for every $d\ge w$. For a prescribed $n\ge w+1$, pick $N=2^m>n$, where $C(N,w)$ follows from \cite[Theorem]{GVB}, and apply the displayed implication successively in degrees $N-1,\ldots,n$. Taking $w=n-1$ proves the conjecture in each degree over $\mathbb C$, since Hypothesis~\ref{Hypothesis:L} is analytic in nature. The relevant Koszul complexes are defined over $\mathbb Q$, whereby Conjecture~\ref{con1} follows over every characteristic-zero field as in Section~\ref{sec:conditional-induction}.

\subsection{Existing results and methods}\label{subsec1.2}
Conjecture~\ref{con1} was originally posed over characteristic $0$ fields $\mathbb{K}$. In general, the conjecture is false over positive characteristic (cf. \cite{GVB} for counterexamples). The authors of \cite{GVB} also show that the truth of Conjecture~\ref{con1} for algebraically closed fields and in a particular degree $d$ is independent of the choice of field, and only depends on the characteristic of the field. The first progress over characteristic $0$ was made for degree $\leq 7$ polynomials via computational methods in \cite{TV}. Soon after, in \cite{GVB} the authors related the conjecture for a general polynomial $P(X)$ of degree $d$ over any algebraically closed field $\mathbb{K}$ to the absence of $\mathbb{K}$-rational points of a weighted projective $\mathbb{Z}$-subscheme $X_d$ of $\mathbb{P}_{\mathbb{Z}}(1, 2, \dots, d-1)$ defined by vanishing of resultants $\Res_X(P, P_i)$ for all $1\leq i<\deg P=d$. Their methods, utilizing reduction modulo prime arguments, successfully proved the conjecture over fields of characteristic $0$ for polynomials of degree $p^k$, $2p^k$ for any prime $p$. Furthermore, by \cite[Proposition 2.2]{GVB}, if Conjecture~\ref{con1} holds for all polynomials of degree $n$ over fields of characteristic $p$, for some $p=0$ or $p$ prime, then Conjecture~\ref{con1} holds true for degree $d$ polynomials over fields of any characteristic except finitely many primes. The finitely many primes $p$ for a given $d$ such that the Conjecture does not hold true in degree $d$ and characteristic $p$ are called \textit{bad primes for $d$}. These methods were reformulated and extended in \cite{DJ} using $p$-adic valuations, where the authors also proved the Conjecture over characteristic $0$ for polynomials of degrees $3 p^k$ and $4p^k$ for primes $p$ greater than $3$ and $4$ respectively (except $p=5, 7$ when $d=4$). 

In this spirit, further results for certain polynomials of degree $5 p^k$ were shown in \cite{SC}. There have also been several computational studies, cf. \cite{CLO}, where the authors verified the conjecture for polynomials of degrees $d=5 p^k,\ 6 p^k$ and $7 p^k$ barring the bad primes in each case, which were also completely classified. These are the only existing positive results towards Conjecture~\ref{con1}. In addition to these, \cite{CLO} also studied obstructions to hypothetical counterexamples and have verified the conjecture for degree $12$. (which is missed by the cases considered above). Computational approaches involving Gröbner bases, even though theoretically possible, get practically infeasible for large degrees due to the complexity of resultants. Alternate approaches to the conjecture have involved analytic tools via the Gauss--Lucas theorem and to a further extent, the theory of Abel-Gontcharoff polynomials: see \cite{Yakubovich14}, \cite{Yakubovich16} and \cite{Massri2023}.

 Recently, Daniel Schaub and Mark Spivakovsky have published a series of papers addressing the conjecture. Their first paper \cite{DM2} provides a partial result towards ``independence" of the higher resultants $\Res_X(P, P_i)$ for $1\leq i\leq d-1$. Their second paper \cite{DM} attempts to extend the reduction modulo prime techniques of existing positive results by obtaining a non-exhaustive list of bad primes $p$ for each $d>1$. The conjecture over characteristic $0$ fields follows for degree $n$ polynomials if and only if there are finitely many bad primes for $p$, which signifies the importance of studying bad primes. However, the bad primes $p$ for a given $d$ can be very large as demonstrated in \cite{GVB}. This has seemed to hinder a purely algebraic approach to the conjecture so far. 

More recently a new approach to Conjecture~\ref{con1} has been initiated by the present author in \cite{SG}. The basis of the main results of \cite{SG} is the following reformulation of Conjecture~\ref{con1} over any algebraically closed field $\mathbb{K}$ and degree $d\geq 3$.

\begin{proposition}[{\cite[Proposition~5.2]{SG}}]\label{Propref}
Let $\mathbb K$ be algebraically closed and let $d\ge3$. Conjecture~\ref{con1} holds in degree $d$ over $\mathbb K$ if and only if, for every $1\le j_1,\ldots,j_{d-1}\le d$, the sequence
\[
 \bigl(\Phi^\#_{d-1,j_i}(HD_{d-1}^{i-1}\mathbf x_{d-1})\mid 1\le i\le d-1\bigr)
\]
is regular in $\mathbb K[x_1,\ldots,x_{d-1}]$.
\end{proposition}
We refer the reader to Section~\ref{sec:notations} for the relevant notation. Using this reformulation, \cite[Theorems~A and~B]{SG} gives an explicit integer $B(d)$ such that $\dim(X_d\times_{\mathbb Z}\mathbb K)\le B(d)-2$ for every algebraically closed field $\mathbb K$ of characteristic zero and $B(d)\le\max\{2,\lfloor(d-1)/2\rfloor\}$, where $\dim\emptyset=-1$. In particular, $X_d(\mathbb K)$ is finite in the degrees for which $B(d)\le2$, and is empty when $B(d)=1$. Furthermore, utilizing \cite[Proposition~5.2]{SG} Schaub and Spivakovsky provide a computable criterion for bad primes for each $d\geq 3$ in their most recent paper \cite{DM3}. This enables them to obtain an upper bound (\cite[Corollary~$3.2$]{DM3}) on the set of bad primes $\mathcal{P}(d)$ for a given $d$ assuming that Conjecture~\ref{con1} holds true in degree $d$ over characteristic $0$.

\subsection{Historical note}\label{subsec1.3}
Besides being of independent interest, Conjecture~\ref{con1} originated in the context of Casas-Alvero's work \cite{CA} on higher-order polar germs of complex analytic plane curve singularities. This problem has also been considered among a list of forty questions on singularities of algebraic varieties (\cite[Problem~$14^{\star\star}$]{HHJS}) compiled by Hauser and Schicho. The following \cite{CA2} is the original context of the conjecture. The main result of \cite{CA} is a factorization theorem for higher-order polars of an irreducible germ of a complex analytic plane curve, generalizing the known one for first order polars. As a corollary, Casas-Alvero then computed the intersection multiplicity of the $r^{th}$ polar and the original germs, which when viewed as a function in $r$ is defined to be the Plücker function of the germ. The Plücker function can be defined even for a non-irreducible germ. The question that led Casas-Alvero to formulate Conjecture~\ref{con1} is that if the Plücker function of a not necessarily irreducible germ is equal to the Plücker function of an irreducible germ, then is the original germ irreducible to start with. If yes, then one would obtain an easy criterion for irreducibility of complex power series in two variables involving the intersection multiplicities of the germ and its polars.

\subsection*{Acknowledgment}
The author would like to thank Mark Spivakovsky for carefully reading through several drafts, suggesting improvements and for providing valuable feedback at various stages. The author also thanks Daniel Schaub and Jacob Tsimerman for comments and questions, which helped identify gaps. The author would also like to thank his advisors Max Lieblich and Farbod Shokrieh for their constant encouragement and support.  This work was partially supported by NSF CAREER DMS-2044564 and NSF FRG DMS-2151718
grants.

\subsection*{AI Statement} The main strategy of this paper is the same as that in the previous versions, which were developed without any LLM assistance. In this paper we remedy a flaw in the previous argument by assuming an explicit analytic hypothesis (Hypothesis~\ref{Hypothesis:L}). This hypothesis was identified and isolated from various analytic and partial results produced by a sustained interaction with AI models over the past few months. In particular, AI assistance was used to obtain the estimates of Section~\ref{sec:evidence}, proofreading the entire paper and formatting.

\section{Global notations and definitions}\label{sec:notations}
\let\thetheorem\oldthetheorem

We follow the notation introduced in \cite{SG}, which we recall here
for convenience.  Unless specified otherwise, $\mathbb K$ is a field
of characteristic zero, and all coordinate variables have degree one.
\begin{enumerate}[(i)]
\item For a univariate polynomial
$f(X)=a_dX^d+a_{d-1}X^{d-1}+\dots+a_0\in\mathbb K[X]$, we denote
its $i$th Hasse--Schmidt derivative by $f_i(X)$, for $0\le i\le d$; thus
\[
 f_i(X)=\binom di a_dX^{d-i}+\binom{d-1}i a_{d-1}X^{d-1-i}
                         +\dots+\binom ii a_i.
\]
Equivalently, $f_i(X)$ is the coefficient of $T^i$ in $f(X+T)$.

\item We denote the $i$th multivariate Hasse--Schmidt operator on
$\mathbb K[x_1,\dots,x_d]$ by $HD_d^i$.  On a monomial it is defined by
\begin{equation}\label{EqnHS1}
 HD_d^i(x_1^{\alpha_1}\dots x_d^{\alpha_d})
 =\sum_{\beta_1+\dots+\beta_d=i}
       \binom{\alpha_1}{\beta_1}\dots\binom{\alpha_d}{\beta_d}
       x_1^{\alpha_1-\beta_1}\dots x_d^{\alpha_d-\beta_d},
\end{equation}
where $0\le \beta_l\le\alpha_l$ in the sum.  We extend this definition
$\mathbb K$-linearly, and put $HD_d^i=0$ when $i<0$.

\item For $d\ge1$ and $1\le j\le d+1$, define the
$\mathbb K$-algebra automorphism $\Phi^\#_{d,j}$ of
$\mathbb K[x_1,\dots,x_d]$ as follows.  Set
$\Phi^\#_{d,d+1}=\operatorname{id}$, and, for $1\le j\le d$, set
\begin{equation}\label{EqnMor}
 \Phi^\#_{d,j}(x_l)=
 \begin{cases}
 x_l-x_j,&l\ne j,\\
 -x_j,&l=j.
 \end{cases}
\end{equation}
The map $\Phi^\#_{d,j}$ is an involution, as one verifies on
the polynomial generators, and hence is an automorphism.
When the polynomial ring is evident from context, we abbreviate
$\Phi^\#_{d,j}$ as $\Phi^\#_j$.

\item Let $\mathbf x_d=x_1x_2\cdots x_d$.  If $e_k$ denotes the
elementary symmetric polynomial of degree $k$, then
\[
 \mathbf x_d=e_d(x_1,\dots,x_d),\qquad
 HD_d^i\mathbf x_d=e_{d-i}(x_1,\dots,x_d).
\]
Write $R_d=\mathbb K[x_1,\ldots,x_d]$.  We use $e_0=1$ and $e_k=0$
for $k<0$ or $k$ larger than the number of entries.
With the convention $x_{d+1}=0$ for an auxiliary zero label in $R_d$,
the polynomials of interest in this paper are
\[
 \Phi^\#_{d,j}(HD_d^{i-1}\mathbf x_d)
   =e_{d+1-i}(x_1-x_j,\dots,x_{d+1}-x_j).
\]
The entry indexed by $j$ is zero and may be omitted.  If $j\le d$,
the last entry is $-x_j$, as in \eqref{EqnMor}.

\item For a choice of indices $1\le j_1,\dots,j_d\le d+1$, write
$\mathbf j=(j_1,\dots,j_d)$ and define
\[
 \begin{split}
 \mathcal S_d(j_1,\dots,j_d)
   &:=\Phi^\#_{d,j_1}(HD_d^{0}\mathbf x_d),\ \Phi^\#_{d,j_2}(HD_d^{1}\mathbf x_d),\dots,\ \Phi^\#_{d,j_d}(HD_d^{d-1}\mathbf x_d)
                                            \\
 \mathcal I_d(j_1,\dots,j_d)
   &:=\bigl(\Phi^\#_{d,j_i}(HD_d^{i-1}\mathbf x_d)
                              \mid1\le i\le d\bigr)
                              \subseteq\mathbb K[x_1,\dots,x_d].
 \end{split}
\]
We abbreviate the ordered sequence and its ideal by
$\mathcal S_d(\mathbf j)$ and $\mathcal I_d(\mathbf j)$.
For a word $\mathbf j=(j_1,\ldots,j_d)$, put
$W_{\mathbf j}=\{j_1,\ldots,j_d\}$ and $|\mathbf j|=|W_{\mathbf j}|$.
Thus $|\mathbf j|$ counts distinct labels, not positions.

\item We denote the Koszul complex of elements $a_1,\ldots,a_r$
in a ring $R$ by $\operatorname K_{R,\bullet}(a_1,\ldots,a_r)$.
\end{enumerate}

We will use the following form of the Koszul criterion in both
standard and weighted polynomial rings, see
\cite[Section~1.6]{BH} and \cite[Lemma~15.31.7]{stacks-project}.
\begin{lemma}\label{lem:koszul-criterion}
Let $R$ be a polynomial ring over a field, with strictly positive
degrees assigned to its variables, and let $a_1,\dots,a_r$ be
homogeneous elements of positive degree.  Then $a_1,\dots,a_r$ is a
regular sequence if and only if
$H_1(\operatorname K_{R,\bullet}(a_1,\dots,a_r))=0$.
In this case all positive Koszul homology vanishes, and every
permutation of the sequence is regular.
\end{lemma}

\section{Elementary identities and the initial cases}
\label{sec:elementary}

We first record a formula for the polynomials
$\Phi^\#_{n,j}(HD_n^{i-1}\mathbf x_n)$, and recall the known
initial case involving at most two selected indices.

\begin{lemma}\label{Lemma:1}
    For any $1\leq i\leq n$ and $1\leq j\leq n+1$, in the ring
    $\mathbb{K}[x_1,\dots, x_n]$, we have
    \begin{equation}\label{Equation:Lemma1}
        \Phi^{\#}_{n,j}(HD^{i-1}_n\mathbf{x}_n)
        =\sum_{l=0}^{n+1-i}\binom{n+1-l}{n+1-i-l}
        (-x_j)^{n+1-i-l}e_l(x_1,\dots, x_n).
    \end{equation}
    Here we use the convention that $x_{n+1}=0$ in
    $\mathbb{K}[x_1,\dots, x_n]$.
\end{lemma}

\begin{proof}
    If $j=n+1$, then $\Phi^{\#}_{n,n+1}$ is the identity map and
    $x_{n+1}=0$, so both sides of \eqref{Equation:Lemma1} equal
    $e_{n+1-i}(x_1,\dots, x_n)$.  We may therefore assume
    $1\leq j\leq n$.

    Note that $\Phi^{\#}_{n,j}(HD^{i-1}_n\mathbf{x}_n)
    =e_{n+1-i}(x_1-x_j,\dots,x_{j-1}-x_j,-x_j,
                              x_{j+1}-x_j,\dots,x_n-x_j)$,
    where $e_{n+1-i}$ is the degree $n+1-i$ elementary symmetric
    polynomial in $n$ variables.  For this proof only, let $z_i:=x_i-x_j$ for all
    $1\leq i\leq n$ with $i\neq j$, and let $z_j=-x_j$.
    Then consider the generating function
    $E_z(T):=\sum_{r=0}^{n}e_r(z_1,\dots,z_n)T^r
             =\prod_{i=1}^{n}(1+z_iT)
             =(1-x_jT)\prod_{i=1,\,i\neq j}^{n}(1+(x_i-x_j)T)
             =(1-x_jT)^n\prod_{i=1,\,i\neq j}^{n}(1+Ux_i)$,
    where $U=T/(1-x_jT)$.  Since
    $(1-x_jT)(1+Ux_j)=1$, we obtain
    $E_z(T)=(1-x_jT)^{n+1}\prod_{i=1}^{n}(1+Ux_i)
          =(1-x_jT)^{n+1}\sum_{i=0}^{n}e_i(x_1,\dots,x_n)U^i
          =\sum_{i=0}^{n}e_i(x_1,\dots,x_n)T^i(1-x_jT)^{n+1-i}$.
    Expanding $(1-x_jT)^{n+1-i}$ using the binomial theorem,
    we obtain
    \begin{align*}
    E_z(T)&=\sum_{r=0}^{n}e_r(z_1,\dots,z_n)T^r
           =\sum_{i=0}^{n}e_i(x_1,\dots,x_n)T^i(1-x_jT)^{n+1-i}\\
          &=\sum_{i=0}^{n}\sum_{k=0}^{n+1-i}
             e_i(x_1,\dots,x_n)T^i
             \binom{n+1-i}{k}(-x_j)^kT^k.
    \end{align*}
    Then comparing coefficients of $T^r$, we see that
    $e_r(z_1,\dots,z_n)
      =\sum_{l=0}^{r}\binom{n+1-l}{r-l}
                       e_l(x_1,\dots,x_n)(-x_j)^{r-l}$.
    Taking $r=n+1-i$, we obtain
    \[
    \Phi^{\#}_{n,j}(HD^{i-1}_n\mathbf{x}_n)
      =e_{n+1-i}(z_1,\dots,z_n)
      =\sum_{l=0}^{n+1-i}\binom{n+1-l}{n+1-i-l}
                       e_l(x_1,\dots,x_n)(-x_j)^{n+1-i-l}.
    \]
\end{proof}

The same generating-function argument gives the general translation
identity
\begin{equation}\label{eq:elementary-translation}
 e_k(t_1-a,\dots,t_N-a)
 =\sum_{l=0}^k\binom{N-l}{k-l}(-a)^{k-l}
                          e_l(t_1,\dots,t_N),\qquad 0\leq k\leq N.
\end{equation}

For the remainder of this section we work in $R_{n-1}$ and use
the auxiliary label $x_n=0$ in its elementary coordinate lists.
For a choice of indices $1\le j_1,\dots,j_{n-1}\le n$, we use the following notation for brevity:
\[
 f_i:=\Phi^\#_{n-1,j_i}(HD_{n-1}^{i-1}\mathbf x_{n-1})
       =e_{n-i}(x_1-x_{j_i},\dots,x_n-x_{j_i})\in R_{n-1}.
\]
Thus $f_1,\dots,f_{n-1}$ is the sequence
$\mathcal S_{n-1}(j_1,\dots,j_{n-1})$ in $R_{n-1}$.  Furthermore, for any $1\leq w\leq n-1$, we will consider the following assertion
\begin{equation}\label{Eq:assertionC(n,w)}
C(n,w): \ \mathcal S_{n-1}(j_1,\dots,j_{n-1})\text{ is a regular sequence in }R_{n-1} \text{ for all } |\mathbf{j}|\leq w \end{equation}
 Since $f_1,\dots, f_{n-1}\in \mathbb{Z}[x_1,\dots, x_{n-1}]$, the assertion $C(n,w)$ is independent of the choice of the characteristic zero field $\mathbb K$. The assertion $C(n,2)$ was established by Draisma and Knopper, see \cite[Proposition~6]{DJ}. Since $\mathcal S_{n-1}(j_1,\dots,j_{n-1})$ is a sequence of homogeneous polynomials in $R_{n-1}$, we observe \begin{equation}\label{Eq:C(n,w)equivalent}
     C(n,w) \text{ is equivalent to } \operatorname{height}\mathcal{I}_{n-1}(\mathbf{j})=n-1 \text{ for all } \mathbf{j}=(j_1,\dots, j_{n-1}) \text{ with } |\mathbf{j}|\leq w,
 \end{equation} 
 where $\mathcal{I}_{n-1}(\mathbf{j})\subseteq R_{n-1}$ is the ideal generated by the sequence $\mathcal S_{n-1}(\mathbf{j})$. The following lemma imposes height bound on $\mathcal{I}_{n-1}(\mathbf{j})\subseteq R_{n-1}$ with $|\mathbf{j}|=w$, under the assumption of $C(n,w-1)$.

\begin{lemma}\label{Lemma:lower-witness-height}
Let $2\le w\le n-1$ and assume $C(n,w-1)$. Then for any choice of indices $\mathbf j=(j_1,\dots, j_{n-1})$ satisfying $|\mathbf{j}|=w$, the ideal $\mathcal{I}_{n-1}(\mathbf{j})=(f_1,\dots,f_{n-1})\subseteq R_{n-1}$ satisfies \[\operatorname{height}\mathcal{I}_{n-1}(\mathbf{j})\ge n-2.\]
\end{lemma}

\begin{proof}
Choose distinct $w_1,w_2\in W_{\mathbf j}$ and set
$\ell=x_{w_2}-x_{w_1}$.  Replace each occurrence of $w_2$ by $w_1$
to obtain a word $\mathbf j'$ using $w-1$ labels.  Modulo $\ell$
each elementary list for $\mathbf j$ agrees with the corresponding list
for $\mathbf j'$.  Thus
\[
 (\mathcal I_{n-1}(\mathbf j),\ell)
 =(\mathcal I_{n-1}(\mathbf j'),\ell).
\]
By $C(n,w-1)$, the quotient by the right-hand ideal is
zero-dimensional.  Since the ideals and $\ell$ are homogeneous,
localization at the homogeneous maximal ideal preserves their
quotient dimensions.  A single element lowers dimension by at most
one \cite[Lemmas~10.117.1 and~10.60.13]{stacks-project}; hence
$\dim R_{n-1}/\mathcal I_{n-1}(\mathbf j)\le1$.
\end{proof}

The following proposition strengthens Lemma~\ref{Lemma:lower-witness-height}
by proving that $f_2,\ldots,f_{n-1}$, together with any difference
between two distinct selected coordinates, is a regular sequence.

\begin{proposition}\label{prop:regularity-with-difference}
Suppose that $2\le w\le n-1$ and that $C(n,w-1)$ holds. Fix a word $\mathbf{j}=(j_1,\dots, j_{n-1})$ with $|\mathbf j|=w$. For
distinct $w_1,w_2\in W_{\mathbf j}=\{j_1,\dots, j_{n-1}\}$, the sequence
\[
 f_2,\ldots,f_{n-1},x_{w_2}-x_{w_1}
\]
is regular in $R_{n-1}$.  In particular, $\mathcal A_{\mathbf j}:=R_{n-1}/(f_2,\dots, f_{n-1})$ is a one-dimensional graded complete intersection.
\end{proposition}

\begin{proof}
Put $\ell=x_{w_2}-x_{w_1}$.  Choose $j'_1=w_1$ and, for $i\ge2$,
replace each occurrence of $w_2$ in $j_i$ by $w_1$, leaving the
other entries unchanged.  The new word uses at most $w-1$ labels.
Writing $f'_i$ for the corresponding generators, $C(n,w-1)$ gives
a regular sequence $f'_2,\ldots,f'_{n-1},f'_1$, by
Lemma~\ref{lem:koszul-criterion}.  Its last member factors as
\[
 f'_1=\ell\prod_{\substack{1\le k\le n\\ k\ne w_1,w_2}}(x_k-x_{w_1}).
\]
A factor of a non-zero divisor is a non-zero divisor, so
$f'_2,\ldots,f'_{n-1},\ell$ is regular as well.  For $i\ge2$ the
elementary lists for $f_i$ and $f'_i$ agree modulo $\ell$.  It
follows that $(f_2,\ldots,f_{n-1},\ell)=(f'_2,\ldots,f'_{n-1},\ell)$
has height $n-1$.  The displayed sequence in the statement is
therefore a homogeneous system of parameters in $R_{n-1}$, hence
regular.  Quotienting by its first $n-2$ elements gives the asserted
one-dimensional complete intersection.
\end{proof}

\section{The coefficient filtration of the Koszul complex}\label{sec:filtration}

\subsection{The two sequences of elementary generators}
\label{subsec2.1}
By Proposition~\ref{Propref}, Conjecture~\ref{con1} in degree $n$
amounts to regularity, for every word $\mathbf j$, of the sequence
\begin{equation}\label{EqProve}
 \mathcal S_{n-1}(\mathbf j)
 =\bigl(\Phi^\#_{n-1,j_i}(HD^{i-1}_{n-1}\mathbf x_{n-1})
                                             \mid 1\le i\le n-1\bigr)
\end{equation}
in $R_{n-1}$ (over an algebraically closed field of characteristic
zero).  We compare it with the first $n-1$ generators of a sequence
in $R_n$.  The filtration argument uses only
\eqref{eq:height-hypothesis} and \eqref{eq:upper-hypothesis} below.

\subsubsection{The hypothesis}\label{hypo}
Fix $n\ge4$ and $3\le w\le n-1$, and assume $C(n+1,w)$ and
$C(n,w-1)$.  To obtain $C(n,w)$, it suffices to consider a word
$\mathbf j$ using exactly $w$ labels in $R_{n-1}$.

\subsubsection{The inductive step}\label{induct}
Write the generators under consideration as
\begin{equation}\label{EqnGoal}
 f_i=\Phi^\#_{n-1,j_i}(HD_{n-1}^{i-1}\mathbf x_{n-1}),
                           \qquad 1\le i\le n-1.
\end{equation}
To adjoin $x_n$ without changing the auxiliary zero label, put
\[
 \widetilde{j}_i=
 \begin{cases}
  j_i,&1\le j_i\le n-1,\\
 n+1,& j_i=n,
 \end{cases}
 \qquad 1\le i\le n-1.
\]
and consider the truncated sequence
\[
 \widehat{\mathcal S}_n(\widetilde{\mathbf{j}}):=\widehat{\mathcal{S}}_n(\widetilde{j}_1,\dots,\widetilde{j}_{n-1})
 :=
 \bigl(
 \Phi^\#_{n,\widetilde{j}_i}(HD^{i-1}_n\mathbf x_n)
 \mid 1\le i\le n-1
 \bigr)
\]
in $R_n$.  Indeed, label $n$ in $R_{n-1}$ and label $n+1$ in $R_n$ both
give the identity automorphism.  The relabeling therefore preserves
the number of distinct labels.  Extend $\widetilde{\mathbf j}$ by
any one of its entries; by $C(n+1,w)$ the resulting length-$n$
sequence is regular, hence so is its prefix
$\widehat{\mathcal S}_n(\widetilde{\mathbf j})$.
Lemma~\ref{Lemma:lower-witness-height} gives
$\operatorname{height}\mathcal I_{n-1}(\mathbf j)\ge n-2$ from $C(n,w-1)$.
We retain the indices $j_i$ for the automorphisms of $R_{n-1}$
and $\widetilde j_i$ for those of $R_n$.  In $R_n$ the convention
$x_{n+1}=0$ is in force.  Since $\widetilde j_i\ne n$,
$x_{\widetilde j_i}$ belongs to the subring $R_{n-1}$; it is
zero when $j_i=n$.  Thus the element denoted $x_{j_i}$ in
the elementary coordinate list of $R_{n-1}$ is represented by
$x_{\widetilde j_i}$ after inclusion in $R_n$.  The indeterminate
$x_n$ of $R_n$ is the newly adjoined variable.

\subsection{The Koszul complexes and their coefficient filtration}
\label{subsec2.2}
For $1\le i\le n-1$,
the elementary identity for adjoining $x_n$ gives
\[
 \Phi^\#_{n,\widetilde j_i}(HD_n^{i-1}\mathbf x_n)
 =(x_n-x_{\widetilde j_i})
        \Phi^\#_{n-1,j_i}(HD_{n-1}^{i-1}\mathbf x_{n-1})
   +\Phi^\#_{n-1,j_i}(HD_{n-1}^{i-2}\mathbf x_{n-1}).
\]
Here $HD_{n-1}^{-1}=0$.  We use the following abbreviations:
\begin{equation}\label{Eqn01}
\begin{split}
F_i= F(i,\widetilde j_i,n)
   &:=\Phi^\#_{n,\widetilde j_i}(HD_n^{i-1}\mathbf x_n)
      =x_nf(i,j_i,n)+g(i,j_i,n),\\
f_i= f(i,j_i,n)&:=\Phi^\#_{n-1,j_i}
                      (HD_{n-1}^{i-1}\mathbf x_{n-1}),\\
g_i=g(i,j_i,n)&:=\Phi^\#_{n-1,j_i}
                      (HD_{n-1}^{i-2}\mathbf x_{n-1})
                       -x_{\widetilde j_i}f(i,j_i,n).
\end{split}
\end{equation}
The coefficients $f(i,j_i,n)$ and $g(i,j_i,n)$ belong to
$R_{n-1}$ and are homogeneous of degrees $n-i$ and $n+1-i$.
Thus $\mathcal S_{n-1}(\mathbf j)$ is the sequence of leading
coefficients of $\widehat{\mathcal S}_n(\widetilde{\mathbf j})$
with respect to $x_n$. The hypotheses used in
this section are
\begin{equation}\label{eq:height-hypothesis}
 \operatorname{height}\mathcal I_{n-1}(\mathbf{j})
 \geq n-2,
\end{equation}
\begin{equation}\label{eq:upper-hypothesis}
 \widehat{\mathcal S}_n(\widetilde{\mathbf j})
 \text{ is a regular sequence in }R_n.
\end{equation}

Consider the Koszul complex of
$\widehat{\mathcal S}_n(\widetilde{\mathbf j})$ in $R_n$:
\begin{equation}\label{Eqn03}
\begin{tikzcd}[column sep=2.5em]
 \dots\arrow[r,"d_{n,4}"] &
 \bigwedge^3 R_n^{\oplus n-1}\arrow[r,"d_{n,3}"] &
 \bigwedge^2 R_n^{\oplus n-1}\arrow[r,"d_{n,2}"] &
 R_n^{\oplus n-1}\arrow[r,"d_{n,1}"] & R_n\arrow[r] &0.
\end{tikzcd}
\end{equation}
By \eqref{eq:upper-hypothesis}, this complex resolves
$R_n/(F_1,\ldots,F_{n-1})$.  We use its first four terms.  All
exact sequences of homology below are taken in the category of
ungraded $R_{n-1}$-modules.

\begin{definition}\label{Def1}
\begin{enumerate}
\item Define
$\widehat{\operatorname{K}}^n(\widetilde{\mathbf j})_{\bullet}$
to be the following chain complex, formed by the terms in
homological degrees $0,1,2,3$ of the Koszul complex~\eqref{Eqn03}
of $\widehat{\mathcal S}_n(\widetilde{\mathbf j})$ in
$R_n=R_{n-1}[x_n]$:
\[
\begin{tikzcd}[column sep=2.2em]
 0\arrow[r] &
 \bigwedge^3R_{n-1}[x_n]^{\oplus n-1}\arrow[r,"d_{n,3}"] &
 \bigwedge^2R_{n-1}[x_n]^{\oplus n-1}\arrow[r,"d_{n,2}"] &
 R_{n-1}[x_n]^{\oplus n-1}\arrow[r,"d_{n,1}"] &
 R_{n-1}[x_n]\arrow[r] &0.
\end{tikzcd}
\]
\item Define
$\operatorname{K}^{n-1}(j_1,\dots,j_{n-1})_{\bullet}$ to be
the following chain complex formed by the terms in homological
degrees $0,1,2,3$ of the Koszul complex of the sequence
$\mathcal S_{n-1}(j_1,\dots,j_{n-1})$ in $R_{n-1}$:
\[
\begin{tikzcd}[column sep=2.5em]
 0\arrow[r] &
 \bigwedge^3R_{n-1}^{\oplus n-1}\arrow[r,"d_{n-1,3}"] &
 \bigwedge^2R_{n-1}^{\oplus n-1}\arrow[r,"d_{n-1,2}"] &
 R_{n-1}^{\oplus n-1}\arrow[r,"d_{n-1,1}"] &
 R_{n-1}\arrow[r] &0.
\end{tikzcd}
\]
\end{enumerate}
In both complexes the rightmost nonzero term has homological
degree zero.  Their homology in degrees $0,1,2$ agrees with that
of the corresponding full Koszul complex.
\end{definition}
We will now consider the truncated complex
$\widehat{\operatorname{K}}^n(\widetilde{\mathbf j})_{\bullet}$ as
a complex of $R_{n-1}$-modules, via the natural embedding
$R_{n-1}\hookrightarrow R_n$, and construct a filtration by
subcomplexes of $R_{n-1}$-modules.

\subsubsection{Filtration of the truncated Koszul complex
$\widehat{\operatorname{K}}^n(\widetilde{\mathbf j})_{\bullet}$}

For any $k\geq0$, let $R_{n-1}[x_n]_k$ be the $R_{n-1}$-submodule
of $R_n$ consisting of polynomials of degree at most $k$ in
$x_n$.  We set $R_{n-1}[x_n]_k=0$ if $k<0$.  Let
$e_1,\dots,e_{n-1}$ be the standard basis of $R_n^{\oplus n-1}$.
For any $p\geq0$ and integer $k$, let $\bigwedge^pR_{n-1}[x_n]_k^{\oplus n-1}$ be the
$R_{n-1}$-submodule of $\bigwedge^pR_n^{\oplus n-1}$ consisting
of sums
\[
 \sum_{1\leq i_1<\dots<i_p\leq n-1}
 h_{i_1\dots i_p}e_{i_1}\wedge\dots\wedge e_{i_p},
 \qquad h_{i_1\dots i_p}\in R_{n-1}[x_n]_k.
\]
Thus the displayed notation stands for
$R_{n-1}[x_n]_k\otimes_{R_{n-1}}\bigwedge^p R_{n-1}^{\oplus n-1}$.
By \eqref{Eqn01}, each $F_i$ is linear in $x_n$.  A Koszul differential
therefore increases the $x_n$-degree of its coefficients by at most one.
This gives the following subcomplexes.
\begin{definition}\label{Def2}
\begin{enumerate}
\item For $k\geq1$, define
$\widehat{\operatorname{K}}^n_k(\widetilde{\mathbf j})_{\bullet}$
to be the following chain complex:
\[
\begin{tikzcd}[column sep=2.1em]
 0\arrow[r] &
 \bigwedge^3R_{n-1}[x_n]_{k-3}^{\oplus n-1}\arrow[r,"d_{n,3}"] &
 \bigwedge^2R_{n-1}[x_n]_{k-2}^{\oplus n-1}\arrow[r,"d_{n,2}"] &
 R_{n-1}[x_n]_{k-1}^{\oplus n-1}\arrow[r,"d_{n,1}"] &
 R_{n-1}[x_n]_k\arrow[r] &0.
\end{tikzcd}
\]
The differentials are the restrictions of those on
$\widehat{\operatorname K}^n(\widetilde{\mathbf j})_\bullet$.
\item Let $M\subseteq R_{n-1}^{\oplus n-1}$ be the image of
the map
$d_{n-1,2}:\bigwedge^2R_{n-1}^{\oplus n-1}
\rightarrow R_{n-1}^{\oplus n-1}$ in
$\operatorname{K}^{n-1}(j_1,\dots,j_{n-1})_{\bullet}$ of
Definition~\ref{Def1}(2).  Define
$\widehat{\operatorname{K}}^n_0(\widetilde{\mathbf j})_{\bullet}$
to be the following complex:
\[
\begin{tikzcd}[column sep=2.5em]
 0\arrow[r] &0\arrow[r] &0\arrow[r] &
 M\arrow[r,"d_{n,1}"] & R_{n-1}[x_n]_0\arrow[r] &0,
\end{tikzcd}
\]
where $M$ lies in homological degree one and $d_{n,1}$ is the
restriction of the map $R_n^{\oplus n-1}\to R_n$
in $\widehat{\operatorname K}^n(\widetilde{\mathbf j})_\bullet$.
It takes values in $R_{n-1}$: for $b=(b_i)\in M$ we have
$\sum_i b_i f_i=0$, and hence $\sum_i b_iF_i=\sum_i b_i g_i$.
\end{enumerate}
\end{definition}
We have natural inclusion chain maps $\iota_{k,\bullet}:
 \widehat{\operatorname{K}}^n_{k-1}(\widetilde{\mathbf j})_{\bullet}
 \hookrightarrow
 \widehat{\operatorname{K}}^n_k(\widetilde{\mathbf j})_{\bullet}$ for all $k\geq1$. The map
$\lambda_{n,k}:R_{n-1}[x_n]_k\rightarrow R_{n-1}$ sending a
polynomial to its coefficient of $x_n^k$ is an $R_{n-1}$-linear
homomorphism for each $k\geq0$.  On a module of exterior coefficients this map acts on each polynomial
coefficient, leaving its wedge basis vector unchanged.  We write the
resulting $R_{n-1}$-linear maps as
\begin{equation}\label{Eqn04}
\begin{aligned}
 \lambda_{n,k}^{\oplus n-1}: R_{n-1}[x_n]_k^{\oplus n-1}
                         \rightarrow R_{n-1}^{\oplus n-1},\quad
 \lambda_{n,k}^{[p]}:
 \bigwedge^pR_{n-1}[x_n]_k^{\oplus n-1}
                         \rightarrow\bigwedge^pR_{n-1}^{\oplus n-1}.
\end{aligned}
\end{equation}
Explicitly, $\lambda_{n,k}^{\oplus n-1}$ sends
$\sum_{i=1}^{n-1}h_ie_i$ to
$\sum_{i=1}^{n-1}\lambda_{n,k}(h_i)e_i$, and we define
\[
 (\lambda_{n,k}^{[p]})
 \left(\sum_{i_1<\dots<i_p}h_{i_1\dots i_p}
                  e_{i_1}\wedge\dots\wedge e_{i_p}\right)
 =\sum_{i_1<\dots<i_p}\lambda_{n,k}(h_{i_1\dots i_p})
                  e_{i_1}\wedge\dots\wedge e_{i_p}.
\]

Since $F_i=x_nf_i+g_i$, taking the coefficient of $x_n^{k-p+1}$
after applying the Koszul differential in degree $p$ amounts to
applying the differential for the $f_i$ to the coefficients of
$x_n^{k-p}$ in degree $p$.  Thus these homomorphisms induce a
surjective chain map $\Lambda_{n,k,\bullet}:\widehat{\operatorname{K}}^n_k(\widetilde{\mathbf j})_{\bullet}
 \rightarrow \operatorname{K}^{n-1}(j_1,\dots,j_{n-1})_{\bullet}$
of $R_{n-1}$-complexes for all $k\geq3$ as follows:
\begin{equation}\label{Eqn05}
\begin{tikzcd}[column sep=2.1em,row sep=3.1em]
 0\arrow[r] &
 \bigwedge^3R_{n-1}[x_n]_{k-3}^{\oplus n-1}
   \arrow[r,"d_{n,3}"]
   \arrow[d,"{\lambda_{n,k-3}^{[3]}}"'] &
 \bigwedge^2R_{n-1}[x_n]_{k-2}^{\oplus n-1}
   \arrow[r,"d_{n,2}"]
   \arrow[d,"{\lambda_{n,k-2}^{[2]}}"'] &
 R_{n-1}[x_n]_{k-1}^{\oplus n-1}
   \arrow[r,"d_{n,1}"]
   \arrow[d,"{\lambda_{n,k-1}^{\oplus n-1}}"'] &
 R_{n-1}[x_n]_k\arrow[r]\arrow[d,"\lambda_{n,k}"'] &0\\
 0\arrow[r] &
 \bigwedge^3R_{n-1}^{\oplus n-1}\arrow[r,"d_{n-1,3}"] &
 \bigwedge^2R_{n-1}^{\oplus n-1}\arrow[r,"d_{n-1,2}"] &
 R_{n-1}^{\oplus n-1}\arrow[r,"d_{n-1,1}"] &
 R_{n-1}\arrow[r] &0.
\end{tikzcd}
\end{equation}
In each degree its kernel consists of the coefficients of one lower
$x_n$-degree.  We obtain the short exact
sequence of $R_{n-1}$-chain complexes for $k\geq3$:
\begin{equation}\label{Eqn06}
\begin{tikzcd}[column sep=2.4em]
 0\arrow[r] &
 \widehat{\operatorname{K}}^n_{k-1}(\widetilde{\mathbf j})_{\bullet}
   \arrow[r,"\iota_{k,\bullet}"] &
 \widehat{\operatorname{K}}^n_k(\widetilde{\mathbf j})_{\bullet}
   \arrow[r,"\Lambda_{n,k,\bullet}"] &
 \operatorname{K}^{n-1}(j_1,\dots,j_{n-1})_{\bullet}
   \arrow[r] &0.
\end{tikzcd}
\end{equation}

For brevity, we will henceforth suppress the indices
$\mathbf j$ and $\widetilde{\mathbf j}$ when referring to the complexes in
Definitions~\ref{Def1} and~\ref{Def2}, when it is clear from
context. There are analogous short exact sequences for
$k=1,2$. Let
$q:R_{n-1}^{\oplus n-1}\rightarrow R_{n-1}^{\oplus n-1}/M$
be the natural quotient map. Since
$M=\operatorname{im}d_{n-1,2}\subseteq\ker d_{n-1,1}$,
the differential $d_{n-1,1}$ induces a homomorphism
$\overline d_{n-1,1}:R_{n-1}^{\oplus n-1}/M\rightarrow R_{n-1}$.
The cokernel of $\iota_{1,\bullet}$ is given by the lower row
of the following diagram:
\begin{equation}\label{EqnCoker1}
\begin{tikzcd}[column sep=2.4em,row sep=3.1em]
 \widehat{\operatorname{K}}^n_{1,\bullet}:
   \arrow[d,"q_{1,\bullet}"'] &
 0\arrow[r] &0\arrow[r] &0\arrow[r] &
 R_{n-1}[x_n]_0^{\oplus n-1}
   \arrow[r,"d_{n,1}"]\arrow[d,"q"'] &
 R_{n-1}[x_n]_1\arrow[r]\arrow[d,"\lambda_{n,1}"'] &0\\
 \operatorname{Coker}\iota_{1,\bullet}: &
 0\arrow[r] &0\arrow[r] &0\arrow[r] &
 R_{n-1}^{\oplus n-1}/M
   \arrow[r,"\overline d_{n-1,1}"] &
 R_{n-1}\arrow[r] &0.
\end{tikzcd}
\end{equation}
Here we identify $R_{n-1}[x_n]_0^{\oplus n-1}$ with
$R_{n-1}^{\oplus n-1}$. Similarly, the cokernel of
$\iota_{2,\bullet}$ is given by the lower row of
\begin{equation}\label{EqnCoker2}
\begin{tikzcd}[column sep=1.6em,row sep=3.1em]
 \widehat{\operatorname{K}}^n_{2,\bullet}:
   \arrow[d,"q_{2,\bullet}"'] &
 0\arrow[r] &0\arrow[r] &
 \bigwedge^2R_{n-1}[x_n]_0^{\oplus n-1}
   \arrow[r,"d_{n,2}"]
   \arrow[d,"{\lambda_{n,0}^{[2]}}"'] &
 R_{n-1}[x_n]_1^{\oplus n-1}
   \arrow[r,"d_{n,1}"]
   \arrow[d,"{\lambda_{n,1}^{\oplus n-1}}"'] &
 R_{n-1}[x_n]_2\arrow[r]\arrow[d,"\lambda_{n,2}"'] &0\\
 \operatorname{Coker}\iota_{2,\bullet}: &
 0\arrow[r] &0\arrow[r] &
 \bigwedge^2R_{n-1}^{\oplus n-1}
   \arrow[r,"d_{n-1,2}"] &
 R_{n-1}^{\oplus n-1}\arrow[r,"d_{n-1,1}"] &
 R_{n-1}\arrow[r] &0.
\end{tikzcd}
\end{equation}
In both diagrams, the rightmost nonzero term has homological
degree zero, and all omitted terms are zero. The diagrams
commute by~\eqref{Eqn01}. Their vertical maps are surjective,
and their kernels are respectively
$\widehat{\operatorname{K}}^n_{0,\bullet}$ and
$\widehat{\operatorname{K}}^n_{1,\bullet}$. Thus, for $k=1,2$,
we obtain the following short exact sequence of
$R_{n-1}$-chain complexes:
\begin{equation}\label{EqnCokerExact}
\begin{tikzcd}[column sep=2.5em]
 0\arrow[r] &
 \widehat{\operatorname{K}}^n_{k-1,\bullet}
   \arrow[r,"\iota_{k,\bullet}"] &
 \widehat{\operatorname{K}}^n_{k,\bullet}
   \arrow[r,"q_{k,\bullet}"] &
 \operatorname{Coker}\iota_{k,\bullet}\arrow[r] &0.
\end{tikzcd}
\end{equation}

Consequently, we have obtained the following filtration of
$\widehat{\operatorname{K}}^n(\widetilde{\mathbf j})_{\bullet}$
by subcomplexes of $R_{n-1}$-modules, along with the cokernels
of successive inclusions.
\begin{equation}\label{Eqn09}
\begin{tikzcd}[column sep=1.7em,row sep=2.8em]
 \widehat{\operatorname{K}}^n_{0,\bullet}
  \arrow[r,hook,"\iota_{1,\bullet}"] &
 \widehat{\operatorname{K}}^n_{1,\bullet}
  \arrow[r,hook,"\iota_{2,\bullet}"]\arrow[d] &
 \widehat{\operatorname{K}}^n_{2,\bullet}
  \arrow[r,hook,"\iota_{3,\bullet}"]\arrow[d] &
 \widehat{\operatorname{K}}^n_{3,\bullet}
  \arrow[r,hook]\arrow[d] &\dots\arrow[r,hook] &
 \widehat{\operatorname{K}}^n_{k,\bullet}
  \arrow[r,hook]\arrow[d] &\dots\arrow[r,hook] &
 \widehat{\operatorname{K}}^n_{\bullet}\\
 &\operatorname{Coker}\iota_{1,\bullet}
 &\operatorname{Coker}\iota_{2,\bullet}
 &\operatorname{K}^{n-1}_{\bullet}
 &&\operatorname{K}^{n-1}_{\bullet}.
\end{tikzcd}
\end{equation}
We can write $\widehat{\operatorname{K}}^n_\bullet$ as the direct
limit $\varinjlim_{l\geq k}\widehat{\operatorname{K}}^n_{l,\bullet}$
for any $k\geq0$.  Filtered direct limits of modules are exact
\cite[Lemma~10.8.8]{stacks-project} and therefore commute with homology,
so for all $k\geq0$ we have
\begin{equation}\label{Eqn13}
 H_{\star}(\widehat{\operatorname{K}}^n_\bullet)
 =\varinjlim_{l\geq k}
             H_{\star}(\widehat{\operatorname{K}}^n_{l,\bullet}).
\end{equation}

\subsubsection{Homology computations}
Consider the two-row matrix
\[
 \mathcal M_n(j_1,\dots,j_{n-1})=
 \begin{bmatrix}
 f(1,j_1,n)&f(2,j_2,n)&\cdots&f(n-1,j_{n-1},n)\\
 g(1,j_1,n)&g(2,j_2,n)&\cdots&g(n-1,j_{n-1},n)
 \end{bmatrix},
\]
and write $D=I_2(\mathcal M_n(j_1,\dots,j_{n-1}))=(f_ig_j-f_jg_i\mid1\le i<j\le n-1)$ for its ideal of $2\times2$ minors.  For $1\le i<j\le n-1$,
\[
 d_{n,1}d_{n-1,2}(e_i\wedge e_j)
       =f_iF_j-f_jF_i=f_ig_j-f_jg_i.
\]
Consequently, $D=d_{n,1}(M)$, where $M$ is the module in
Definition~\ref{Def2}(2). Thus, we obtain
\begin{equation}\label{EqnH0-i}
 \begin{split}
 H_0(\widehat{\operatorname K}^{n}_{\bullet})
   =R_n/(\widehat{\mathcal S}_n(\widetilde{\mathbf j})),\quad
H_0(\widehat{\operatorname K}^{n}_{k,\bullet})=\begin{cases}
R_{n-1}/D\ & k=0,\\
    R_{n-1}[x_n]_k\big/
        d_{n,1}(R_{n-1}[x_n]_{k-1}^{\oplus n-1})\ &k\geq 1
\end{cases}
\end{split}
\end{equation}

\begin{lemma}\label{Lemma:filteredHone}
Under \eqref{eq:height-hypothesis} and \eqref{eq:upper-hypothesis}, we have $H_2(\operatorname K^{n-1}_{\bullet})=0$ and $H_1(\widehat{\operatorname K}^{n}_{k,\bullet})=0$ for any $k\ge2$.
\end{lemma}

\begin{proof}
Since $R_{n-1}$ is Cohen--Macaulay, hypothesis
\eqref{eq:height-hypothesis} gives
\[
 \operatorname{grade}\mathcal I_{n-1}(j_1,\dots,j_{n-1})
 =\operatorname{height}\mathcal I_{n-1}(j_1,\dots,j_{n-1})
 \ge n-2.
\]
By depth sensitivity of Koszul homology
\cite[Theorem~1.6.17]{BH}, the homology of the full Koszul complex of $f_1,\ldots,f_{n-1}$ vanishes in degrees greater than
$(n-1)-(n-2)=1$.  In particular,
$H_2(\operatorname K^{n-1}_{\bullet})=0$, since the truncation
in Definition~\ref{Def1} preserves second homology. For $k\ge3$, the long exact sequence in homology associated to
\eqref{Eqn06} contains
\[
 H_2(\operatorname K^{n-1}_{\bullet})
 \rightarrow H_1(\widehat{\operatorname K}^{n}_{k-1,\bullet})
 \xrightarrow{\iota_{k,\star}}
 H_1(\widehat{\operatorname K}^{n}_{k,\bullet}),
\]
whose first term is zero.  Thus the transition maps on first homology are injective from
level $2$ onwards.  By \eqref{Eqn13} their direct limit is
$H_1(\widehat{\operatorname K}^n_\bullet)=0$, the last equality
following from \eqref{eq:upper-hypothesis}.  Every term from level
$2$ onwards must vanish.
\end{proof}

We next compute the kernel of each map on zeroth homology of the filtration.

\begin{theorem}\label{MainTheorem}\label{Prop:injection}
Under \eqref{eq:height-hypothesis} and \eqref{eq:upper-hypothesis},
for every $k\ge1$ there is an exact sequence of
$R_{n-1}$-modules
\begin{equation}\label{eq:intro-exact}
\begin{split}
 0\rightarrow H_1(\operatorname K^{n-1}_{\bullet})
 \xrightarrow{\delta_k}H_0(\widehat{\operatorname K}^{n}_{k-1,\bullet})
 \xrightarrow{\iota_{k,\star}}H_0(\widehat{\operatorname K}^{n}_{k,\bullet})\rightarrow
 R_{n-1}/\mathcal I_{n-1}(j_1,\dots,j_{n-1})
 \rightarrow0.
\end{split}
\end{equation}
In particular, $\iota_{1,\star}$ is injective if and only if
$\mathcal S_{n-1}(\mathbf j)$ is regular.
\end{theorem}

\begin{proof}
For $k\ge3$, the long exact sequence in homology associated to
\eqref{Eqn06}, together with Lemma~\ref{Lemma:filteredHone},
gives~\eqref{eq:intro-exact}.
For $k=1,2$, the descriptions in \eqref{EqnCoker1} and
\eqref{EqnCoker2} give
\[
\begin{aligned}
 H_1(\operatorname{Coker}\iota_{k,\bullet})
   \simeq H_1(\operatorname K^{n-1}_{\bullet}),\quad
 H_0(\operatorname{Coker}\iota_{k,\bullet})
   \simeq R_{n-1}/\mathcal I_{n-1}(j_1,\dots,j_{n-1}).
\end{aligned}
\]
For $k=2$, Lemma~\ref{Lemma:filteredHone} gives
$H_1(\widehat{\operatorname K}^{n}_{2,\bullet})=0$.
The long exact sequence in homology associated to
\eqref{EqnCokerExact} therefore gives~\eqref{eq:intro-exact}. For $k=1$, it remains to show that the map
\[
 q_{1,\star}:H_1(\widehat{\operatorname K}^{n}_{1,\bullet})
 \rightarrow H_1(\operatorname{Coker}\iota_{1,\bullet})
\]
induced by $q_{1,\bullet}$ in \eqref{EqnCokerExact} is zero. The descriptions
\eqref{EqnCoker1} and \eqref{EqnCoker2} give a chain map
\[
 \pi_\bullet:\operatorname{Coker}\iota_{2,\bullet}
       \rightarrow\operatorname{Coker}\iota_{1,\bullet},
\]
whose components in homological degrees two, one and zero are,
respectively, zero, $q$ and $\operatorname{id}_{R_{n-1}}$.
This is a chain map because $M=\operatorname{im}d_{n-1,2}$
and $\overline d_{n-1,1}\circ q=d_{n-1,1}$. Furthermore,
multiplication by $x_n$ defines a chain map
$\widehat{\operatorname K}^{n}_{1,\bullet}\rightarrow
 \widehat{\operatorname K}^{n}_{2,\bullet}$, and the definitions
of the coefficient maps yield the following commutative diagram:
\[
\begin{tikzcd}[column sep=3.6em,row sep=2.8em]
 \widehat{\operatorname K}^{n}_{1,\bullet}
   \arrow[r,"\times x_n"]\arrow[d,"q_{1,\bullet}"'] &
 \widehat{\operatorname K}^{n}_{2,\bullet}
   \arrow[d,"q_{2,\bullet}"]\\
 \operatorname{Coker}\iota_{1,\bullet} &
 \operatorname{Coker}\iota_{2,\bullet}
   \arrow[l,"\pi_\bullet"']
\end{tikzcd}
\]
Thus $q_{1,\star}$ factors through
$H_1(\widehat{\operatorname K}^{n}_{2,\bullet})$, which is zero
by Lemma~\ref{Lemma:filteredHone}. Applying the long exact
sequence in homology to \eqref{EqnCokerExact} now proves
\eqref{eq:intro-exact} also for $k=1$.The final assertion follows from Lemma~\ref{lem:koszul-criterion}.
\end{proof}

The same criterion holds at every positive level of the filtration.

\begin{corollary}\label{Corollary:filteredcriterion}
Under the hypotheses of Theorem~\ref{MainTheorem}, the
following assertions are equivalent:
\begin{enumerate}[(i)]
\item The sequence $\mathcal S_{n-1}(j_1,\dots,j_{n-1})$
is regular in $R_{n-1}$.
\item The map $\iota_{1,\star}:H_0(\widehat{\operatorname K}^{n}_{0,\bullet})\to H_0(\widehat{\operatorname K}^{n}_{1,\bullet})$ is injective.
\item The map $\iota_{k,\star}:H_0(\widehat{\operatorname K}^{n}_{k-1,\bullet})\to H_0(\widehat{\operatorname K}^{n}_{k,\bullet})$ is injective for some, or equivalently all $k\ge1$.
\end{enumerate}
\end{corollary}

\begin{proof}
By Theorem~\ref{MainTheorem}, the kernel of each
$\iota_{k,\star}:H_0(\widehat{\operatorname K}^{n}_{k-1,\bullet})\to H_0(\widehat{\operatorname K}^{n}_{k,\bullet})$ is isomorphic to $H_1(\operatorname K^{n-1}_{\bullet})$.  By
Lemma~\ref{lem:koszul-criterion}, its vanishing is equivalent to
regularity of $\mathcal S_{n-1}(j_1,\dots,j_{n-1})$.  This proves the equivalence of the three assertions.
\end{proof}

\subsection{The first Koszul homology}\label{subsec:first-koszul-homology}
Fix $2\le w\le n-1$ and a word $\mathbf j$ using exactly $w$ labels.
Under $C(n,w-1)$, Proposition~\ref{prop:regularity-with-difference}
makes $\mathcal A_{\mathbf j}=R_{n-1}/(f_2,\ldots,f_{n-1})$ a
one-dimensional graded complete intersection.  We identify the first
Koszul homology with the annihilator of the product of the unused
factors of $f_1$.  Here
$U_{\mathbf j}=\{1,\ldots,n\}\setminus W_{\mathbf j}$ is the set of
unused labels.

We fix a coordinate convention for the rest of the paper.  The inclusion
$R_{n-1}=\mathbb K[x_1,\ldots,x_{n-1}]\hookrightarrow
R_n=R_{n-1}[x_n]$ preserves every element of $R_{n-1}$ and adjoins
the independent variable $x_n$.  The original labels remain
$1,\ldots,n$: their coordinate list in $R_{n-1}$ is
$(x_1,\ldots,x_{n-1},0)$.  In $R_n$ this same list is indexed by
$1,\ldots,n-1,n+1$, with $x_{n+1}=0$, as in Section~\ref{induct}.
In particular, $x_{\widetilde j_1}$ is an element of $R_{n-1}$:
it equals $x_{j_1}$ if $j_1<n$, and equals $0$ if $j_1=n$. For each original label $1\le k\le n$, define the following fixed
element of $R_{n-1}$, and use the same symbol for its image in $R_n$:
\begin{equation}\label{eq:translated-coordinates}
 z_k:=
 \begin{cases}
  x_k-x_{\widetilde j_1},&1\le k<n,\\
  -x_{\widetilde j_1},&k=n.
 \end{cases}
\end{equation}
Thus $z_{j_1}=0$, $z_n=-x_{\widetilde j_1}$, and
\[
 z_{j_i}=x_{\widetilde j_i}-x_{\widetilde j_1}
 \qquad(1\le i\le n-1).
\]
The subscript of every $z_k$ always denotes an original label in $\{1,\ldots,n\}$, irrespective of the ambient ring being $R_{n-1}$ or $R_n$.  In particular, $z_n$ is independent of $x_n$, when considered in $R_n$.  If $j_1=n$, the convention gives $z_k=x_k$ for $k<n$ and $z_n=0$. The translated coordinate corresponding to the new variable is $T:=x_n-x_{\widetilde j_1}=x_n+z_n$. The inverse relations are
\[
 x_k=z_k-z_n\quad(1\le k<n),\qquad x_n=T-z_n.
\]
Consequently, the $n-1$ elements $z_k$ with $k\ne j_1$ are
polynomial coordinates on $R_{n-1}$, and adjoining $T$ gives
polynomial coordinates on $R_n$:
\begin{equation}\label{Eq:newcoordinates}
 R_{n-1}=\mathbb K[(z_k)_{k\ne j_1}],\qquad
 R_n=R_{n-1}[T]=\mathbb K[(z_k)_{k\ne j_1},T].
\end{equation}
In the order of the labels $1,\ldots,n+1$ of $R_n$, the translated
coordinate list is $(z_1,\ldots,z_{n-1},T,z_n)$.  Elementary
symmetric polynomials are unchanged if this list is reordered as
$(z_1,\ldots,z_n,T)$.  These conventions apply throughout the
remaining sections. Let $\ell_k=z_k$ for $k\ne j_1$.  Thus
\begin{equation}\label{eq:elementary-generators-z}
 f_i=e_{n-i}(z_1-z_{j_i}, \dots, z_n-z_{j_i}),
 \ \text{ for all }1\leq i\leq n-1,\quad \text{in particular, }\
 f_1=\prod_{k\in\{1,\ldots,n\}\setminus\{j_1\}}\ell_k.
\end{equation}
Let $q_U:=\prod_{u\in U_{\mathbf j}}\ell_u$, the product of the
linear forms corresponding to the unused indices.

\begin{proposition}\label{prop:koszul-annihilator}
Assuming $C(n,w-1)$, for any fixed word $\mathbf{j}:=(j_1,\dots, j_{n-1})$ with $|\mathbf{j}|=w$, there is an isomorphism of $R_{n-1}$-modules
\[
 H_1(\operatorname K^{n-1}_\bullet)
 \simeq\operatorname{Ann}_{\mathcal A_{\mathbf j}}(q_U).
\]
Moreover, if one assumes $C(n+1,w)$ then for every $k\ge1$,
\begin{equation}\label{eq:filtered-annihilator}
 \operatorname{Ann}_{\mathcal A_{\mathbf j}}(q_U)\simeq \ker\left(\iota_{k,\star}:H_0(\widehat{\operatorname K}^{n}_{k-1,\bullet})
 \rightarrow H_0(\widehat{\operatorname K}^{n}_{k,\bullet})\right).
\end{equation}
\end{proposition}

\begin{proof}
Permuting the generators does not change Koszul homology.
The full Koszul complex on $f_2,\ldots,f_{n-1},f_1$ is the
mapping cone of multiplication by $f_1$ on the Koszul complex
of $f_2,\ldots,f_{n-1}$ by \cite[Lemma~15.29.8]{stacks-project}.
The latter complex resolves $\mathcal A_{\mathbf j}=R_{n-1}/(f_2,\dots, f_{n-1})$ by
Proposition~\ref{prop:regularity-with-difference}. The homology sequence of the cone and Definition~\ref{Def1} give the exact sequence
\[
 0\rightarrow H_1(\operatorname K^{n-1}_\bullet)
 \rightarrow\mathcal A_{\mathbf j}
 \xrightarrow{\ f_1\ }\mathcal A_{\mathbf j}
 \rightarrow R_{n-1}/(f_1,\ldots,f_{n-1})
 \rightarrow0.
\]
By Proposition~\ref{prop:regularity-with-difference}, multiplication
by each selected factor $\ell_k$, $k\in W_{\mathbf j}\setminus\{j_1\}$,
is injective on $\mathcal A_{\mathbf j}$.  Consequently
\[
 \operatorname{Ann}_{\mathcal A_{\mathbf j}}(f_1)
 =\operatorname{Ann}_{\mathcal A_{\mathbf j}}(q_U)
\]
by \eqref{eq:elementary-generators-z}.  If also $C(n+1,w)$ holds,
the hypotheses of Theorem~\ref{MainTheorem} follow from
Section~\ref{induct}, and \eqref{eq:intro-exact} identifies this
annihilator with $\ker\iota_{k,\star}$ for every $k\ge1$.
\end{proof}

It suffices to test one unused factor:
\begin{lemma}\label{lem:unused-permutations}
$\operatorname{Ann}_{\mathcal A_{\mathbf j}}(q_U)=0$ if and only if $\operatorname{Ann}_{\mathcal A_{\mathbf j}}(\ell_u)=0$ for some $u\in U_\mathbf{j}$.
\end{lemma}

\begin{proof}
Any permutation $\tau$ of $U_{\mathbf j}$, extended by the identity on
$W_{\mathbf j}$, induces the automorphism $z_k\mapsto z_{\tau(k)}$ of
$R_{n-1}$.  It fixes every $f_i$ by \eqref{eq:elementary-generators-z}
and hence descends to $\mathcal A_{\mathbf j}$, where it conjugates
multiplication by $\ell_u$ to multiplication by $\ell_{\tau(u)}$.
Thus one unused factor is regular if and only if all are.  Their
product $q_U$ is regular precisely in this case.
\end{proof}

The preceding lemma also relates the obstruction to Lefschetz properties of graded complete intersections.  A linear form is a \emph{weak Lefschetz element} if multiplication by it has maximal rank between every pair of
consecutive graded components and a \emph{strong Lefschetz element} if the same holds for all its positive powers, see
\cite[Definition~2.7]{PalezzatoTorielli2020Lefschetz}
for these definitions without an Artinian assumption.
Under $C(n,w-1)$, any $\ell_k$ with
$k\in W_{\mathbf j}\setminus\{j_1\}$ multiplies injectively on
$\mathcal A_{\mathbf j}$.  Its Hilbert function is therefore
nondecreasing.  Thus an unused factor
$\ell_u$ is regular on $\mathcal A_{\mathbf j}$
if and only if it is a weak Lefschetz element, or equivalently a strong
Lefschetz element, of $\mathcal A_{\mathbf j}$.

There is also an equivalent formulation in terms of the usual Artinian
complete intersections considered in
\cite[Section~1]{HarimaMiglioreNagelWatanabe2003Lefschetz}.
Choose $w_0\in W_{\mathbf j}\setminus\{j_1\}$ and set
\[
 s=\frac{(n-2)(n-3)}2,\qquad
 \mathcal B_{\mathbf j,w_0}
   :=\mathcal A_{\mathbf j}/(\ell_{w_0}^{\,s+2}).
\]
By Proposition~\ref{prop:regularity-with-difference}, $\ell_{w_0}$ is
regular on $\mathcal A_{\mathbf j}$, so $\mathcal B_{\mathbf j,w_0}$
is an Artinian complete intersection of socle degree $2s+1$:
its generator degrees are $1,\ldots,n-2,s+2$, so their degrees
minus one sum to $s+(s+1)$.  For every $u\in U_{\mathbf j}$ we have
the following criterion.

\noindent\textbf{Claim}: $\operatorname{Ann}_{\mathcal A_{\mathbf j}}(\ell_u)=0\quad\Longleftrightarrow\quad\text{the image of $\ell_u$ is a weak Lefschetz element of
 $\mathcal B_{\mathbf j,w_0}$.}$
\begin{proof}[Proof of the claim]
The numerator of the Hilbert series
\[
 \operatorname{Hilb}_{\mathcal A_{\mathbf j}}(t)
 =\frac{\prod_{r=1}^{n-2}(1+t+\cdots+t^{r-1})}{1-t}
\]
has degree $s$. Hence
$\dim_{\mathbb K}(\mathcal A_{\mathbf j})_d
 =\dim_{\mathbb K}(\mathcal A_{\mathbf j})_s$ for $d\ge s$.
Since $\ell_{w_0}$ is regular, it follows that
\[
 \times\ell_{w_0}^{\,d-s}:
 (\mathcal A_{\mathbf j})_s
 \xrightarrow{\ \sim\ }(\mathcal A_{\mathbf j})_d
 \qquad(d\ge s).
\]
Moreover, the quotient map identifies
$(\mathcal B_{\mathbf j,w_0})_d$ with
$(\mathcal A_{\mathbf j})_d$ for $0\le d\le s+1$.

Suppose that the image of $\ell_u$ is a weak Lefschetz element of
$\mathcal B_{\mathbf j,w_0}$. Since the dimensions in degrees $s$
and $s+1$ agree, multiplication
\[
 \ell_u:(\mathcal A_{\mathbf j})_s
 \rightarrow(\mathcal A_{\mathbf j})_{s+1}
\]
is an isomorphism. Let $a\in(\mathcal A_{\mathbf j})_d$ satisfy
$\ell_u a=0$. If $d\le s$, then
$\ell_u(\ell_{w_0}^{s-d}a)=0$, so $\ell_{w_0}^{s-d}a=0$ and hence $a=0$.
If $d>s$, write $a=\ell_{w_0}^{d-s}a'$ with
$a'\in(\mathcal A_{\mathbf j})_s$. Then
$\ell_{w_0}^{d-s}\ell_u a'=0$, whence $\ell_u a'=0$ and $a'=0$.
Thus $\operatorname{Ann}_{\mathcal A_{\mathbf j}}(\ell_u)=0$.

Conversely, suppose that $\ell_u$ is regular on
$\mathcal A_{\mathbf j}$. Its multiplication maps on
$\mathcal A_{\mathbf j}$ are injective in every degree and,
by the dimension equality above, bijective in degrees $d\ge s$.
The induced maps on $\mathcal B_{\mathbf j,w_0}$ are therefore
injective for $d\le s$, since the relevant graded components are
unchanged, and surjective for $d\ge s$, since surjectivity
descends to the quotient. Hence the image of $\ell_u$ is a weak
Lefschetz element of $\mathcal B_{\mathbf j,w_0}$.
\end{proof}

In particular, the image of $\ell_{w_0}$ is a weak Lefschetz element of
$\mathcal B_{\mathbf j,w_0}$.  Indeed, the converse argument uses only regularity on
$\mathcal A_{\mathbf j}$ and applies to $\ell_{w_0}$ as well.
The question is whether the specified forms $\ell_u$ also have this
property.  It is therefore
a question about the non-Lefschetz locus studied in
\cite[Section~2]{BoijMiglioreMiroRoigNagel2018Locus}, namely the locus
where the middle multiplication matrix has zero determinant.
The permutations in Lemma~\ref{lem:unused-permutations} fix $\ell_{w_0}$,
hence descend to $\mathcal B_{\mathbf j,w_0}$ and conjugate these
matrices for the unused indices.  Thus this Artinian formulation
preserves both the symmetry of Lemma~\ref{lem:unused-permutations}
and the Koszul obstruction in~\eqref{eq:filtered-annihilator}.

\section{Uniform local estimates and the homogeneous obstruction}\label{sec:analytic-estimate}

We work over $\mathbb C$.  Fix $n\ge3$, a word
$\mathbf j=(j_1,\ldots,j_{n-1})\in\{1,\ldots,n\}^{n-1}$ and an unused
index $u\in U_{\mathbf j}$.  In the coordinates
\eqref{eq:translated-coordinates}, $z_{j_1}=0$ identically.
On the hyperplane $z_u=0$ write
\begin{equation}\label{Eqn02}
 q=n-2,\qquad
 \mathbf{y}=(z_k)_{k\in\{1,\ldots,n\}\setminus\{j_1,u\}},\qquad
 d_i=n-i\quad(1\le i\le n-1).
\end{equation}
Recall from \eqref{eq:elementary-generators-z} that
$f_i=e_{d_i}((z_k-z_{j_i})_{k=1}^n)$ for $1\le i\le n-1$.
For $\mathbf c=(c_2,\ldots,c_{n-1})\in\mathbb C^q$, put
\begin{equation}\label{eq:coefficient-family}
 \begin{split}
 \gamma_i&=e_{d_i+1}\bigl((z_k-z_{j_i})_{k=1}^n\bigr)-z_{j_i}f_i
           =g_i-z_nf_i\quad(1\le i\le n-1),\\
 Q_{i,u,\mathbf c}&=(f_i+c_i\gamma_i)|_{z_u=0}
           =\bar f_i+c_i\bar\gamma_i\quad(2\le i\le n-1).
 \end{split}
\end{equation}
Here $\bar f_i=f_i|_{z_u=0}$ and $\bar\gamma_i=\gamma_i|_{z_u=0}$.
We write $Q_{u,\mathbf c}=(Q_{2,u,\mathbf c},\ldots,Q_{n-1,u,\mathbf c})$:
the index $u$ specifies the hyperplane of restriction, and $i$
specifies the component.  If $c_i\ne0$, one also has the identity
\[
 Q_{i,u,\mathbf c}
 =c_i e_{d_i+1}\bigl((z_k-z_{j_i})_{k=1}^n,
                         c_i^{-1}-z_{j_i}\bigr)|_{z_u=0}.
\]
Define the finite set
\begin{equation}\label{eq:analytic-coefficient-grid}
 \mathscr G_n
 =\prod_{i=2}^{n-1}
       \left\{\frac{n+1-i}{r}:1\le r\le i-1\right\}
 \subset\mathbb Q_{>0}^q.
\end{equation}
It has $q!$ elements. Norms and distances are Euclidean unless otherwise specified; $\|\mathbf{y}\|_\infty=\max_{1\le k\le q}|y_k|$.

To verify the relation with $g_i$ in \eqref{eq:coefficient-family},
recall that $T=x_n+z_n$ is the additional coordinate in $R_n$.
Adjoining the entry $T-z_{j_i}$ to the elementary symmetric polynomial gives
\[
 F_i=e_{d_i+1}\bigl((z_k-z_{j_i})_{k=1}^n,T-z_{j_i}\bigr)
     =e_{d_i+1}\bigl((z_k-z_{j_i})_{k=1}^n\bigr)+(T-z_{j_i})f_i.
\]
By \eqref{Eqn01}, $F_i=(T-z_n)f_i+g_i$.  Comparison of constant
terms in $T$ yields $\gamma_i=g_i-z_nf_i=g_i+x_{\widetilde j_1}f_i$.
For $i=1$ the list contains $z_{j_1}=0$, so $\gamma_1=0$ and
\begin{equation}\label{eq:adjoined-coordinate}
 F_i=Tf_i+\gamma_i\quad(1\le i\le n-1),\qquad
 F_1=T\prod_{k\ne j_1}z_k .
\end{equation}
The coordinates $((z_k)_{k\ne j_1},T)$ identify
$R_n/(T-1,z_u)$ with $\mathbb C[\mathbf{y}]$.  In this quotient
$F_i=\bar f_i+\bar\gamma_i=Q_{i,u,\mathbf1}$ for $2\le i\le n-1$.
In the original coordinates the specialization $T=1$ is $x_n=1-z_n$.

For $n\ge4$, Hypothesis~\ref{Hypothesis:L} requires that, whenever $0$ is an
isolated zero of $Q_{u,\mathbf1}$, the estimate
\eqref{eq:finite-grid-estimate} hold uniformly for $\mathbf c$ in a
complex neighborhood of $\mathscr G_n$.  Its constants may depend
on $n,\mathbf j,u$.  Since the grid is finite, one may equivalently
allow separate constants and exponents at its centers: for $\|\mathbf{y}\|<1$,
take the largest of the exponents and the smallest of the lower-bound
constants and radii.

\subsection*{The homogeneous obstruction in the selected coordinates}

We first express the obstruction of Section~\ref{subsec:first-koszul-homology}
in the selected coordinates.  Work in $R_{n-1}$ before imposing $z_u=0$.
For a word using $w=|W_{\mathbf j}|\ge2$ labels, write
\[
 \mathbf{v}=(z_k)_{k\in W_{\mathbf j}\setminus\{j_1\}},\qquad
 \mathbf{u}=(z_k)_{k\in U_{\mathbf j}},\qquad
 s=w-1,\quad r=n-w,\quad q=s+r-1.
\]
Here $\mathbf{u}$ denotes the tuple of unused coordinates, whereas
$u$ is the fixed unused label.  Thus $R_{n-1}=\mathbb C[\mathbf{v},\mathbf{u}]$.
Put $a_d=z_{j_{n-d}}$ for
$1\le d\le q$.  Following the elementary recursions in
\cite[Lemma~5.19 and (5.22)]{SG}, set $D_0=U_0=1$ and
\begin{equation}\label{eq:gw-general-D-U-recursion}
 D_d=-\sum_{e=0}^{d-1}\binom{n-e}{d-e}(-a_d)^{d-e}D_e,\qquad
 U_d=D_d-\sum_{e=1}^{d}e_e(\mathbf{v})U_{d-e}.
\end{equation}
These are homogeneous polynomials of degree $d$ in $\mathbf{v}$.
Writing $E_d=e_d(\mathbf{v},\mathbf{u})$, translation by $-a_d$ and convolution
of elementary symmetric functions give, respectively,
\[
 \begin{aligned}
 f_{n-d}&=\sum_{e=1}^{d}\binom{n-e}{d-e}(-a_d)^{d-e}(E_e-D_e),\\
 E_d-D_d&
   =\sum_{e=1}^{d}e_{d-e}(\mathbf{v})\bigl(e_e(\mathbf{u})-U_e\bigr).
 \end{aligned}
\]
Here $e_e(\mathbf{u})=0$ for $e>r$.  Both changes of generators are triangular
with diagonal entries $1$.  Consequently, setting
$B_{\mathbf j}:=\mathbb C[\mathbf{v}]/(U_{r+1},\ldots,U_q)$
(with no relations when $r=q$), we obtain
\begin{equation}\label{eq:gw-general-homogeneous-normal-form}
 \mathcal A_{\mathbf j}
 =R_{n-1}/(f_2,\ldots,f_{n-1})
 \simeq B_{\mathbf j}[\mathbf{u}]/(e_1(\mathbf{u})-U_1,\ldots,e_r(\mathbf{u})-U_r).
\end{equation}
This is free of rank $r!$ over $B_{\mathbf j}$ by the splitting
algebra basis theorem \cite[Proposition~1.10]{EkedahlLaksov2005}.
In this presentation $q_U=\prod_{k\in U_{\mathbf j}}z_k=e_r(\mathbf{u})=U_r$.
Multiplication by $q_U$ therefore acts as multiplication by $U_r$ on
each of the $r!$ free summands, giving an isomorphism of ungraded
$B_{\mathbf j}$-modules
\begin{equation}\label{eq:gw-arbitrary-word-obstruction-reduction}
 \operatorname{Ann}_{\mathcal A_{\mathbf j}}(q_U)
 \simeq \operatorname{Ann}_{B_{\mathbf j}}(U_r)^{\oplus r!}.
\end{equation}
If $f_2,\ldots,f_{n-1}$ is a regular sequence and the selected
factors $z_k$ ($k\in W_{\mathbf j}\setminus\{j_1\}$) are regular on
$\mathcal A_{\mathbf j}$, this is also the full Koszul obstruction:
\[
 H_1\bigl(\operatorname K_{R_{n-1},\bullet}(f_1,\ldots,f_{n-1})\bigr)
 \simeq \operatorname{Ann}_{B_{\mathbf j}}(U_r)^{\oplus r!}
 \qquad\text{as ungraded $B_{\mathbf j}$-modules}.
\]
Indeed, the Koszul complex on $f_2,\ldots,f_{n-1}$ resolves
$\mathcal A_{\mathbf j}$, so adjoining $f_1$ gives
$H_1\simeq\operatorname{Ann}_{\mathcal A_{\mathbf j}}(f_1)$.
The selected factors of $f_1$ act injectively, leaving precisely
the annihilator of $q_U$.  These hypotheses hold under $C(n,w-1)$ by
Proposition~\ref{prop:regularity-with-difference}.

On $z_u=0$ only $r-1$ unused coordinates remain, so $e_r(\mathbf{u})=0$.
Thus the homogeneous equations project onto
$U_r=\cdots=U_q=0$ in the $s$ selected coordinates: every solution
lifts by taking the remaining $r-1$ unused coordinates with elementary
symmetric functions $U_1,\ldots,U_{r-1}$.  Over $\mathbf v=0$ all
$U_d$ of positive degree vanish, so every unused coordinate also
vanishes.  Thus the restricted forms have no nonzero common zero
if and only if $U_r,\ldots,U_q$ have no nonzero common zero in
the selected coordinates.

\subsection{The finite coefficient set and the exponent}

We return to an arbitrary word $\mathbf j$ and unused index $u$.
The following identity gives the coefficient set \eqref{eq:analytic-coefficient-grid}.

\begin{lemma}\label{lem:analytic-mobius}
Put $\varphi_t(z)=(1-t)z/(1-tz)$ and
$D(t)=\prod_{k=1}^n(1-tz_k)$.  For $j\in\{1,\ldots,n\}$, $\alpha=z_j$ and
$0\le d\le n$, one has, wherever the denominators are nonzero,
\begin{equation}\label{eq:analytic-mobius}
\begin{split}
 &e_d\bigl(\varphi_t(z_k)-\varphi_t(\alpha):1\le k\le n\bigr)\\
 &\quad=\frac{(1-t)^d}{(1-t\alpha)^dD(t)}
   \sum_{r=0}^{n-d}(-t)^r\binom{d+r}{d}
      (1-t\alpha)^{n-d-r}
      e_{d+r}\bigl(z_k-\alpha:1\le k\le n\bigr).
\end{split}
\end{equation}
\end{lemma}
\begin{proof}
Set $[n]=\{1,\ldots,n\}$ and $\xi_k=z_k-\alpha$.  Since
\[
 \varphi_t(z_k)-\varphi_t(\alpha)
 =\frac{(1-t)\xi_k}{(1-tz_k)(1-t\alpha)},
\]
the defining sum over $d$-element subsets gives
\[
 e_d\bigl(\varphi_t(z_k)-\varphi_t(\alpha):k\in[n]\bigr)
 =\frac{(1-t)^d}{(1-t\alpha)^dD(t)}
   \sum_{\substack{S\subseteq[n]\\ |S|=d}}
       \prod_{k\in S}\xi_k\prod_{k\notin S}(1-tz_k).
\]
For a fixed $S$, each factor indexed by its complement equals
$(1-t\alpha)-t\xi_k$.  Choose $-t\xi_k$ for the indices in a subset
$R\subseteq[n]\setminus S$ and $1-t\alpha$ for the others.  If $|R|=r$,
the corresponding term, including the factors indexed by $S$, is
\[
 (-t)^r(1-t\alpha)^{n-d-r}\prod_{k\in S\cup R}\xi_k .
\]
For fixed $r$, each $(d+r)$-element set $U$ occurs as $S\cup R$
exactly $\binom{d+r}{d}$ times: choose $S\subseteq U$ of size $d$, and then
$R=U\setminus S$ is determined.  The sum of the products over all such
$U$ is $e_{d+r}(\xi_1,\ldots,\xi_n)$.  Summing over $0\le r\le n-d$
gives \eqref{eq:analytic-mobius}, including $d=0,n$ by the
empty-product convention.
\end{proof}

\begin{theorem}\label{thm:analytic-equivalence}
For a fixed pair $(\mathbf j,u)$ with $u\in U_\mathbf{j}$, the following are equivalent.
\begin{enumerate}[(i)]
\item The forms
$f_2|_{z_u=0},\ldots,f_{n-1}|_{z_u=0}$ have only the origin
as a common zero in $\mathbb C^q$.
\item There are $A,B,\varepsilon>0$ for which
\eqref{eq:finite-grid-estimate} holds.
\item There is an $\varepsilon>0$ such that, whenever
$\|\mathbf{y}\|<\varepsilon$ and $\operatorname{dist}(\mathbf c,\mathscr G_n)<\varepsilon$,
the equality $Q_{u,\mathbf c}(\mathbf{y})=0$ implies $\mathbf{y}=0$.
\item For every compact $K\subset\mathbb C^q$, there are
$A_K,R_K>0$ such that
\begin{equation}\label{eq:compact-coefficient-estimate}
 \|Q_{u,\mathbf c}(\mathbf{y})\|\ge A_K\|\mathbf{y}\|^q
 \qquad (\mathbf c\in K,\ \|\mathbf{y}\|<R_K).
\end{equation}
\end{enumerate}
\end{theorem}
\begin{proof}
Assume \textup{(i)} and let $K\subset\mathbb C^q$ be nonempty and compact.
Set
\[
 m_0=\min_{\|\omega\|=1}\|(\bar f_i(\omega))_{i=2}^{n-1}\|>0,\qquad
 M_K=\max_{\substack{\mathbf c\in K\\ \|\omega\|=1}}
          \|(c_i\bar\gamma_i(\omega))_{i=2}^{n-1}\|.
\]
Choose $0<R_K\le1$ with $R_KM_K\le m_0/2$.
For $0<r<R_K$ the $i$th component of $Q_{u,\mathbf c}(r\omega)$
equals $r^{d_i}(\bar f_i(\omega)+rc_i\bar\gamma_i(\omega))$.
Since $r^{d_i}\ge r^q$, the triangle inequality gives
\[
 \|Q_{u,\mathbf c}(r\omega)\|
 \ge r^q\|(\bar f_i(\omega)+rc_i\bar\gamma_i(\omega))_{i=2}^{n-1}\|
 \ge r^q(m_0-rM_K)\ge (m_0/2)r^q .
\]
This proves \textup{(iv)}.  A compact neighborhood of the finite set
$\mathscr G_n$ gives \textup{(iv)}$\Rightarrow$\textup{(ii)}, after
shrinking $\varepsilon$; \textup{(ii)}$\Rightarrow$\textup{(iii)} follows
from the positivity of the right-hand side for $\mathbf{y}\ne0$.

We prove \textup{(iii)}$\rightarrow$\textup{(i)} by
contradiction.  Let $0\ne p\in\mathbb C^q$ be a common zero
of $\bar f_2,\ldots,\bar f_{n-1}$.  By homogeneity we may rescale $p$
so that $z_{k_0}(p)=1$ for some $k_0\notin\{j_1,u\}$.  Put
$\zeta_k=z_k(p)$ for all $k$, with $\zeta_{j_1}=\zeta_u=0$.
Define $\eta(t)\in\mathbb C^q$ near $t=0$ by
\[
 z_k(\eta(t))=\varphi_t(\zeta_k)
 \qquad\bigl(k\notin\{j_1,u\}\bigr).
\]
Since $\varphi_t(0)=0$ and $\varphi_t(1)=1$ near $t=0$,
the two omitted labelled values remain zero, and
$z_{k_0}(\eta(t))=1$.  Set
\[
 D_p(t)=\prod_{k=1}^n(1-t\zeta_k).
\]
In particular, $D_p(0)=1$. Fix $i\in\{2,\ldots,n-1\}$, and abbreviate
\[
 d=d_i=n-i,\qquad
 \alpha=\zeta_{j_i},\qquad
 A_r=e_{d+r}(\zeta_k-\alpha:1\le k\le n)
 \quad(0\le r\le i).
\]
Here $A_0=\bar f_i(p)=0$, while $A_i=0$ because the entry
indexed by $j_i$ is zero. Suppose first that some $A_r$ is nonzero, and let $r_i=\min\{r\in\{1,\ldots,i-1\}:A_r\ne0\}$. Define
\[
 \begin{aligned}
 H_i(t)&=\sum_{r=r_i}^{i-1}(-t)^{r-r_i}
     \binom{d+r}{d}(1-t\alpha)^{i-r}A_r,\\
 J_i(t)&=\sum_{r=r_i}^{i-1}(-t)^{r-r_i}
     \binom{d+r}{d+1}(1-t\alpha)^{i-r}A_r .
 \end{aligned}
\]
Applying Lemma~\ref{lem:analytic-mobius} in degree $d$
gives the exact identity
\[
 \bar f_i(\eta(t))
 =
 (-t)^{r_i}
 \frac{(1-t)^d}{(1-t\alpha)^dD_p(t)}H_i(t).
\]
Applying the same lemma in degree $d+1$ gives
\[
 e_{d+1}\bigl(
 z_k(\eta(t))-z_{j_i}(\eta(t)):1\le k\le n
 \bigr)
 =
 (-t)^{r_i-1}
 \frac{(1-t)^{d+1}}
      {(1-t\alpha)^{d+1}D_p(t)}J_i(t).
\]
Since $z_{j_i}(\eta(t))=\varphi_t(\alpha)=(1-t)\alpha/(1-t\alpha)$, the definition of $\gamma_i$ therefore yields
\[
 \bar\gamma_i(\eta(t))
 =
 (-t)^{r_i-1}
 \frac{(1-t)^{d+1}}
      {(1-t\alpha)^{d+1}D_p(t)}
 \bigl(J_i(t)+\alpha tH_i(t)\bigr).
\]

The constant terms satisfy $H_i(0)=\binom{d+r_i}{d}A_{r_i}\ne0$ and $J_i(0)=\binom{d+r_i}{d+1}A_{r_i}\ne0$. Thus,
\[
 c_i(t)=
 \frac{1-t\alpha}{1-t}\,
 \frac{H_i(t)}{J_i(t)+\alpha tH_i(t)}
\]
is holomorphic near $t=0$, with $c_i(0)= \binom{d+r_i}{d}/\binom{d+r_i}{d+1} =\frac{d_i+1}{r_i}$. The preceding exact identities give
\begin{align*}
 \bar f_i(\eta(t))+t c_i(t)\bar\gamma_i(\eta(t))=
 \frac{(1-t)^dH_i(t)}
      {(1-t\alpha)^dD_p(t)}
 \left((-t)^{r_i}+t(-t)^{r_i-1}\right)
 =0.
\end{align*}

If $A_r=0$ for every $0\le r\le i$, the same
applications of Lemma~\ref{lem:analytic-mobius} show that
$\bar f_i(\eta(t))$ and $\bar\gamma_i(\eta(t))$ both vanish
identically.  In this case set $c_i(t)=d_i+1$.
Thus in both cases
$\bar f_i(\eta(t))+t c_i(t)\bar\gamma_i(\eta(t))=0$.  Shrinking to a common neighborhood of $t=0$
gives a holomorphic $\mathbf c(t)=(c_2(t),\ldots,c_{n-1}(t))$ with
$\mathbf c(0)\in\mathscr G_n$.  By homogeneity,
\[
 Q_{i,u,\mathbf c(t)}(t\eta(t))
 =t^{d_i}\bigl(\bar f_i(\eta(t))+tc_i(t)\bar\gamma_i(\eta(t))\bigr)=0
 \qquad (2\le i\le n-1).
\]
But $z_{k_0}(t\eta(t))=t$.  Thus these nonzero zeros tend to the
origin while $\mathbf c(t)\to\mathbf c(0)\in\mathscr G_n$, contrary
to \textup{(iii)}.
\end{proof}

\begin{remark}\label{rem:analytic-quantifiers}
Although condition (iv) of Theorem~\ref{thm:analytic-equivalence}
appears stronger than (ii), the two conditions are equivalent for
the family \eqref{eq:coefficient-family}.
The proof even permits $A_K=m_0/2$ independently of $K$; the
radius $R_K$ must be allowed to depend on $K$.
For example, when $n=3$, take $z_{j_1}=z_u=0$ and let the
remaining coordinate and the second selected coordinate be $y$.
Writing $c=c_2$, the family has just one component:
\[
 Q_{u,c}(y)=Q_{2,u,c}(y)=-2y+3cy^2.
\]
For $c\ne0$ its additional zero is $2/(3c)$, which tends to
zero as $c\to+\infty$.  Consequently no fixed
positive radius works for all compact coefficient sets.
\end{remark}

Under $C(n,w-1)$, the first change of selected index reduces the set of
centers in Theorem~\ref{thm:analytic-equivalence}.  Here $w=|W_{\mathbf j}|$.  Let $\nu$ be the least degree
$2\le\nu\le q$ with $j_{n-\nu}\ne j_{n-1}$, if such a degree exists
(as it does when $|W_{\mathbf j}|\ge3$).  Put
\[
 a=z_{j_{n-1}},\qquad b=z_{j_{n-\nu}},\qquad i_*=n-\nu+1.
\]
Here $i_*$ indexes the last generator with selected index $j_{n-1}$
before the first change: $f_{i_*}$ has degree $\nu-1$ and center $a$,
whereas $f_{i_*-1}$ has degree $\nu$ and center $b$.
The next proposition forces $c_{i_*}(0)=\nu$ for the arc constructed
in the proof of \textup{(iii)}$\Rightarrow$\textup{(i)} of
Theorem~\ref{thm:analytic-equivalence}.  Only centers with $c_{i_*}=\nu$ need to be tested;
homogeneity allows the test to be made within the same hyperplane.

\begin{proposition}\label{prop:first-change-coefficient}
In $R_{n-1}$ one has
\begin{equation}\label{eq:first-change-coefficient}
\begin{split}
 \gamma_{i_*}=f_{n-\nu}-\binom n\nu(a-b)^\nu-a f_{i_*}-\sum_{e=1}^{\nu-1}\binom{n-e}{\nu-e}(a-b)^{\nu-e}f_{n-e}.
\end{split}
\end{equation}
Consequently, $\gamma_{i_*}=-\binom n\nu(a-b)^\nu$ in
$\mathcal A_{\mathbf j}$.  If $|W_{\mathbf j}|=w$ and $C(n,w-1)$
holds, then $\gamma_{i_*}$ is a non-zero divisor in this ring.
\end{proposition}
\begin{proof}
Write $E_e=e_e((z_k-a)_{k=1}^n)$.  Translating by $a-b$ in
\eqref{eq:elementary-translation} gives
\[
 f_{n-\nu}=\sum_{e=0}^{\nu}\binom{n-e}{\nu-e}(a-b)^{\nu-e}E_e.
\]
Here $E_0=1$, $E_e=f_{n-e}$ for $1\le e<\nu$, and
$E_\nu=\gamma_{i_*}+a f_{i_*}$, giving \eqref{eq:first-change-coefficient}.
Modulo $(f_2,\ldots,f_{n-1})$ only $-\binom n\nu(a-b)^\nu$
remains.  Under $C(n,w-1)$ the difference $a-b$ is regular by
Proposition~\ref{prop:regularity-with-difference}, and $\binom n\nu\ne0$.
\end{proof}

Now let $\mathcal H_{\mathbf j}=\{\mathbf c\in\mathbb C^q:c_{i_*}=\nu\}$ be a distinguished hyperplane in $\mathbb{C}^q$ and let $\mathscr G^*_{n,\mathbf j}=\mathscr G_n\cap\mathcal H_{\mathbf j}$. Thus, only $r=1$ remains out of the $n-\nu$ choices $\nu/r$ for $c_{i_*}$ in \eqref{eq:analytic-coefficient-grid}, whereby the number of centers drops from $q!$ in $\mathscr G_n$ to $q!/(n-\nu)$ in $\mathscr G^*_{n,\mathbf j}$.

\begin{theorem}\label{thm:first-change-sliced-criterion}
Assume $C(n,w-1)$ and $|W_{\mathbf j}|=w\ge3$.  For each fixed
unused index $u\in U_{\mathbf{j}}$, the following are equivalent:
\begin{enumerate}[(i)]
\item The only common zero of $f_2|_{z_u=0},\ldots,f_{n-1}|_{z_u=0}$ is the origin.
\item There are $A,B,\varepsilon>0$ such that
\begin{equation}\label{eq:first-change-sliced-estimate}
 \|Q_{u,\mathbf c}(\mathbf{y})\|\ge A\|\mathbf{y}\|^B
\end{equation}
whenever $\|\mathbf{y}\|<\varepsilon$, $\mathbf c\in\mathcal H_{\mathbf j}$
and $\operatorname{dist}(\mathbf c,\mathscr G^*_{n,\mathbf j})<\varepsilon$.
\item For some $\varepsilon>0$, every zero $Q_{u,\mathbf c}(\mathbf{y})=0$
in the range specified in \textup{(ii)} satisfies $\mathbf{y}=0$.
\end{enumerate}
When these conditions hold, one may take $B=q$ in \textup{(ii)}.
\end{theorem}
\begin{proof}
Applying Theorem~\ref{thm:analytic-equivalence}\textup{(iv)} on a
compact neighborhood of $\mathscr G^*_{n,\mathbf j}$, gives
\textup{(i)}$\Rightarrow$\textup{(ii)} with $B=q$;
\textup{(ii)}$\Rightarrow$\textup{(iii)} is immediate.

If \textup{(i)} fails, rescale a nonzero common zero $p$ so that one
coordinate equals $1$.  By Proposition~\ref{prop:regularity-with-difference},
the full homogeneous regular sequence $f_2,\ldots,f_{n-1},a-b$ vanishes
only at the origin.  Thus $a(p)\ne b(p)$, and
\eqref{eq:first-change-coefficient} gives $\bar\gamma_{i_*}(p)\ne0$. To $p$, we apply the M\"obius transformation used in the proof of
\textup{(iii)}$\Rightarrow$\textup{(i)} of Theorem~\ref{thm:analytic-equivalence}.
With $\zeta_k=z_k(p)$ and $z_k(\eta(t))=\varphi_t(\zeta_k)$, that
construction gives $\mathbf{y}(t)=t\eta(t)$ and a holomorphic coefficient vector
$\mathbf c(t)$ such that $Q_{u,\mathbf c(t)}(\mathbf{y}(t))=0$ and
$\mathbf c(0)\in\mathscr G_n$.

Consider the component $Q_{i_*,u,\mathbf c}=\bar f_{i_*}+c_{i_*}\bar\gamma_{i_*}$.
For $i=i_*$ one has $d_i=\nu-1$ and $\zeta_{j_i}=a(p)$.
Thus the coefficient $A_1$ defined in that construction is
\[
 A_1=e_\nu(\zeta_k-a(p):1\le k\le n)
    =\bar\gamma_{i_*}(p)+a(p)\bar f_{i_*}(p)
    =\bar\gamma_{i_*}(p)\ne0,
\]
since $\bar f_{i_*}(p)=0$.  Hence the least positive $r$ with $A_r\ne0$
is $r_{i_*}=1$.  The formula for the limiting coefficient now gives
\[
 c_{i_*}(0)=\frac{d_{i_*}+1}{r_{i_*}}=\nu,\qquad
 \mathbf c(0)\in\mathscr G^*_{n,\mathbf j}.
\]

Set $\lambda(t)=c_{i_*}(t)/\nu$, so $\lambda(0)=1$. By homogeneity, we have $Q_{i,u,\mathbf c/\lambda}(\lambda \mathbf{y})=\lambda^{d_i}Q_{i,u,\mathbf c}(\mathbf{y})$ for all $\lambda\ne0$. Thus $\mathbf c(t)/\lambda(t)$ lies in $\mathcal H_{\mathbf j}$ and tends
to the same center, while $\lambda(t)\mathbf{y}(t)$ is a zero tending to the origin.
Its normalized coordinate is $\lambda(t)t\ne0$ for small $t\ne0$,
contradicting \textup{(iii)}.
\end{proof}

\section{The conditional induction}\label{sec:conditional-induction}

We first verify the isolatedness assumption in Hypothesis~\ref{Hypothesis:L}.
Fix $n\ge3$ and a word $\mathbf j$; retain \eqref{Eqn01} and work over
$\mathbb C$.

\begin{lemma}\label{lem:distinguished-isolated}
Suppose that $F_1,\ldots,F_{n-1}$ is a regular sequence in $R_n$.
For every unused index $u\in U_{\mathbf j}$, the origin $\mathbf{0}\in \mathbb{C}^q$ is an
isolated zero of $Q_{u,\mathbf1}$.
\end{lemma}

\begin{proof}
Use the linear coordinates $(T,z)$ on $\mathbb C^n$, where
$z=(z_k)_{k\ne j_1}$, and set $X=V(F_1,\ldots,F_{n-1})$.
The $n-1$ forms are a regular sequence in the $n$-dimensional
polynomial ring $R_n$, so $\dim X=n-(n-1)=1$.

First consider the slice $S=X\cap\{T=1\}$.  Since the $F_i$ are
homogeneous, simultaneous scaling of $T$ and $z$ gives an isomorphism
\[
 \mathbb C^\times\times S \rightarrow X\cap\{T\ne0\},
 \qquad (\lambda,(1,z))\longmapsto(\lambda,\lambda z),
\]
with inverse $(T,z)\mapsto(T,(1,z/T))$.  The right-hand side has
dimension at most one, so $S$ has no positive-dimensional component.
In particular, $S$ has only finitely many points. Now fix an unused index $u$.  On $T=1$, $z_u=0$, the equation
$F_1=0$ holds automatically: by \eqref{eq:adjoined-coordinate},
$F_1=T\prod_{k\ne j_1}z_k$ and $u\ne j_1$.  For $2\le i\le n-1$,
the same identity $F_i=Tf_i+\gamma_i$ gives
\[
 F_i\big|_{T=1,\,z_u=0}
   =\bar f_i+\bar\gamma_i=Q_{i,u,\mathbf1}.
\]
Thus, in the remaining coordinates $\mathbf{y}=(z_k)_{k\notin\{j_1,u\}}$,
the zero set of $Q_{u,\mathbf1}$ is exactly $S\cap\{z_u=0\}$.
It is finite and contains $\mathbf{y}=\mathbf{0}$, since $\bar f_i=f_i|_{z_u=0}$ and
$\bar\gamma_i=\gamma_i|_{z_u=0}$ have positive degree.  Hence $\mathbf{y}=\mathbf{0}$ is an isolated zero of $Q_{u,\mathbf1}$.
\end{proof}

\begin{theorem}\label{thm:conditional-step}
Assume Hypothesis~\ref{Hypothesis:L}. For $n\ge4$ and
$3\le w\le n-1$, one has
\[
 C(n+1,w)\ \text{and}\ C(n,w-1)
 \rightarrow C(n,w).
\]
\end{theorem}

\begin{proof}
It suffices to work over $\mathbb C$. Indeed, the polynomials in the sequence $\mathcal{S}_{n-1}(j_1\dots, j_{n-1})\in R_{n-1}$ are defined over $\mathbb Q$, and flat base change for the associated Koszul complex $K_{\mathbb Q,\bullet}$ yields
\[
 H_1(K_{\mathbb Q,\bullet}\otimes_{\mathbb Q}\mathbb K)
 \simeq H_1(K_{\mathbb Q,\bullet})\otimes_{\mathbb Q}\mathbb K.
\]
every field extension $\mathbb K/\mathbb Q$. By faithful flatness and Lemma~\ref{lem:koszul-criterion},
the assertions in \eqref{Eq:assertionC(n,w)} are therefore
independent of the characteristic-zero field. Assume $C(n+1,w)$ and $C(n,w-1)$ and fix a word $\mathbf j=(j_1,\dots, j_{n-1})$ with
$|\mathbf j|=w$ since regularity of $\mathcal{S}_{n-1}(j_1\dots, j_{n-1})$ with $|\mathbf{j}|<w$ are covered by $C(n,w-1)$. Extend the lifted word $\widetilde{\mathbf j}$ of
Section~\ref{induct} by taking $\widetilde j_n=\widetilde j_1$.
The extended word still uses $w$ indices, so $C(n+1,w)$ gives
regularity of the associated sequence $(F_1,\ldots,F_n)$ and hence of its prefix $\widehat{\mathcal S}_n(\widetilde{\mathbf j})=(F_1,\ldots,F_{n-1})$ in $R_n$.
Moreover, $\operatorname{height}\mathcal I_{n-1}(\mathbf j)\ge n-2$ by Lemma~\ref{Lemma:lower-witness-height}. Thus the hypotheses of Theorem~\ref{MainTheorem} hold. At $k=1$, the exact sequence \eqref{eq:intro-exact} and Proposition~\ref{prop:koszul-annihilator}
identify
\[
 \ker\!\left(
 H_0(\widehat{\operatorname K}^{n}_{0,\bullet})
 \xrightarrow{\iota_{1,\star}}
 H_0(\widehat{\operatorname K}^{n}_{1,\bullet})
 \right)
 \simeq H_1(\operatorname K^{n-1}_\bullet)
 \simeq \operatorname{Ann}_{\mathcal A_{\mathbf j}}(q_U).
\]

Choose $u\in U_{\mathbf j}$ and let $S_u=R_{n-1}/(\ell_u)=\mathbb C[\mathbf{y}]$, where $\ell_u=z_u$ and
$\mathbf{y}$ is as in \eqref{Eqn02}. By Lemma~\ref{lem:distinguished-isolated}, $\mathbf{0}\in \mathbb{C}^q$ is an isolated zero of $Q_{u,\mathbf1}$. The estimate \eqref{eq:finite-grid-estimate} given by Hypothesis~\ref{Hypothesis:L} and implication \textup{(ii)}$\Rightarrow$\textup{(i)} of Theorem~\ref{thm:analytic-equivalence}, we obtain $V(\bar f_2,\ldots,\bar f_{n-1})=\{\mathbf{0}\}$ in $\mathbb C^{n-2}$. The $n-2$ forms $\bar f_2,\ldots,\bar f_{n-1}$ are therefore a homogeneous system of parameters in the polynomial ring $S_u$, and hence a regular sequence.

Now let $\operatorname K_{R_{n-1},\bullet}
(f_2,\ldots,f_{n-1})$ be the Koszul complex of the seuqnece $f_2,\dots, f_{n-1}$ in $R_{n-1}$, which resolves $\mathcal A_{\mathbf j}$ by Proposition~\ref{prop:regularity-with-difference},. Since $\ell_u$ is a coordinate of $R_{n-1}$, multiplication by $\ell_u$ gives the following short exact sequence of complexes
\[
 0\rightarrow \operatorname K_{R_{n-1},\bullet}
(f_2,\ldots,f_{n-1})
 \xrightarrow{\ \ell_u\ } \operatorname K_{R_{n-1},\bullet}
(f_2,\ldots,f_{n-1})
 \rightarrow
 \operatorname K_{S_u,\bullet}(\bar f_2,\ldots,\bar f_{n-1})
 \rightarrow0.
\]
The associated long exact sequence in homology, identifies
\[
 \operatorname{Ann}_{\mathcal A_{\mathbf j}}(\ell_u)
 \simeq H_1\!\left(
   \operatorname K_{S_u,\bullet}(\bar f_2,\ldots,\bar f_{n-1})
 \right)=0,
\]
where the vanishing follows from $\bar f_2,\ldots,\bar f_{n-1}$ being a homogeneous system of parameters in the polynomial ring $S_u$. By Lemma~\ref{lem:unused-permutations}, we therefore obtain
$\operatorname{Ann}_{\mathcal A_{\mathbf j}}(q_U)=0$, whereby $\iota_{1,\star}$ is injective, and
Corollary~\ref{Corollary:filteredcriterion} yields the regularity of
$\mathcal S_{n-1}(\mathbf j)$, as required.
\end{proof}

Thus, we can now prove Conjecture~\ref{con1}, conditional on Hypothesis~\ref{Hypothesis:L}.

\begin{proof}[Proof of Theorem~\ref{MainConditional}]
Assume Hypothesis~\ref{Hypothesis:L} and work first over $\mathbb C$.
We prove by induction on $w\ge2$ that $C(d,w)$ holds for every $d\ge w+1$.
For $w=2$ this is the two-witness result of Draisma and Knopper
\cite[Proposition~6]{DJ}; if all the derivatives have
the same witness, Taylor expansion at that witness gives a pure power.

Let $w\ge3$ and assume that $C(d,w-1)$ holds for all $d\ge w$.
Fix $n\ge w+1$ and choose a power of two $N>n$.
The prime-power theorem of Graf von Bothmer, Labs, Schicho and
van de Woestijne \cite{GVB}, together with Proposition~\ref{Propref},
gives $C(N,N-1)$, hence $C(N,w)$.  Descend from $N$ to $n$.
At each step $d=N-1,\ldots,n$, we have $C(d+1,w)$ from the previous
degree and $C(d,w-1)$ from induction on $w$.  Since $3\le w\le d-1$,
Theorem~\ref{thm:conditional-step} gives $C(d,w)$.  This proves $C(n,w)$, and $n$ was arbitrary. Taking $w=n-1$ yields $C(n,n-1)$, which is equivalent to Conjecture~\ref{con1} in degree $n$, for all $n\ge3$. 
\end{proof}

\section{Uniform estimates for arbitrary words}\label{sec:evidence}

Fix $n\ge3$, a word $\mathbf j$, and an unused index $u\in U_{\mathbf j}$.
Retain $q=n-2$ and $Q_{i,u,\mathbf c}=\bar f_i+c_i\bar\gamma_i$ from
\eqref{eq:coefficient-family}. We first restrict the centers at which
the estimate in Hypothesis~\ref{Hypothesis:L} must be tested.
We then prove an unconditional uniform estimate on a common nonempty
Zariski open set of coefficients and determine which prescribed centers
can belong to such an open set. Nonemptiness of the open set alone
does not imply that it contains a prescribed center.

\subsection{Uniform inequalities for arbitrary words}

Recall from \eqref{eq:gw-general-D-U-recursion}
that the selected coordinate in degree $h$ is $a_h=z_{j_{n-h}}$
for $1\le h\le q$. Set
\[
 \mathscr B_{\mathbf j}
 =\{1\le h<q:j_{n-h-1}\notin\{j_{n-h},\ldots,j_{n-1}\}\}.
\]
Thus $h\in\mathscr B_{\mathbf j}$ exactly when the index selected in
degree $h+1$ is new.  In this case we adjoin the degree-$(h+1)$ form
$\bar\gamma_{n-h}$ from the preceding component.  Consequently, the number of additional forms is $|\{j_2,\ldots,j_{n-1}\}|-1\le |W_{\mathbf j}|-1$.

\begin{theorem}\label{thm:gw-completion}
Let $H_n=(n+3)4^n$ and $\delta_n=(2q4^{n-1}H_n^{2(q-1)})^{-1}$.
For every $\mathbf{y}\ne0$, one has
\begin{equation}\label{eq:gw-completion}
 \max\left\{\max_{1\le h\le q}\frac{|\bar f_{n-h}(\mathbf{y})|}{\|\mathbf{y}\|_\infty^h},
       \ \max_{h\in\mathscr B_{\mathbf j}}
           \frac{|\bar\gamma_{n-h}(\mathbf{y})|}{\|\mathbf{y}\|_\infty^{h+1}}\right\}
 \ge\delta_n .
\end{equation}
The maximum over an empty index set is defined to be zero.  In particular,
for $0<\|\mathbf{y}\|_\infty\le1$,
\begin{equation}\label{eq:gw-completion-unweighted}
 \max\left\{\max_{2\le i\le n-1}|\bar f_i(\mathbf{y})|,\
       \max_{h\in\mathscr B_{\mathbf j}}|\bar\gamma_{n-h}(\mathbf{y})|\right\}
 \ge\delta_n\|\mathbf{y}\|_\infty^q .
\end{equation}
\end{theorem}
\begin{proof}
By homogeneity it suffices to take $\|\mathbf{y}\|_\infty=1$.  Let $\epsilon$
denote the left side of \eqref{eq:gw-completion} at this point.  Since
$z_{j_1}=z_u=0$, we have $\max_k|z_k|=1$ by \eqref{Eqn02}.  As in
Proposition~\ref{prop:first-change-coefficient}, let $a=a_1$ and
$E_d=e_d(z_k-a:1\le k\le n)$, with $E_0=1$.  We prove, by induction on $h$, that
\begin{equation}\label{eq:gw-completion-induction}
 \max\{|E_d|,\ |a_d-a|^h\}\le H_n^{2(h-1)}\epsilon
 \qquad(1\le d\le h).
\end{equation}
For $h=1$, $E_1=\bar f_{n-1}$ and $a_1=a$.  Suppose
\eqref{eq:gw-completion-induction} holds for $h<q$ and choose
$\alpha=a_{h+1}$ if $h\notin\mathscr B_{\mathbf j}$ and $\alpha=a_h$
otherwise.  In either case $\alpha$ is among $a_1,\ldots,a_h$ and
\[
 e_{h+1}(z_k-\alpha:1\le k\le n)
 =\begin{cases}
 \bar f_{n-h-1},&h\notin\mathscr B_{\mathbf j},\\
 \bar\gamma_{n-h}+a_h\bar f_{n-h},&h\in\mathscr B_{\mathbf j}.
 \end{cases}
\]
Its modulus is at most $2\epsilon$ since $|a_h|\le1$.  By
\eqref{eq:elementary-translation},
\[
 e_{h+1}(z_k-\alpha:1\le k\le n)
 =E_{h+1}+\sum_{d=1}^h\binom{n-d}{h+1-d}(a-\alpha)^{h+1-d}E_d
          +\binom n{h+1}(a-\alpha)^{h+1}.
\]
Since $|a-\alpha|\le2$, each coefficient of $E_d$ in the sum is
bounded by $2^n2^{h+1-d}\le4^n$.  Since $\alpha$ was already selected,
the induction hypothesis gives $|a-\alpha|^{h+1}\le2H_n^{2(h-1)}\epsilon$.
It follows that
\[
 \begin{aligned}
 |E_{h+1}|&\le 2\epsilon+4^n\sum_{d=1}^{h}|E_d|
          +\binom n{h+1}|a-\alpha|^{h+1}\\
 &\le(2+h4^n+2^{n+1})H_n^{2(h-1)}\epsilon
 \le H_n^{2h-1}\epsilon .
 \end{aligned}
\]
For $1\le d\le h$, we have $|a_d-a|^{h+1}\le2H_n^{2(h-1)}\epsilon$.
For $d=h+1$, we apply the same translation formula at $a_{h+1}$,
where the left side equals $\bar f_{n-h-1}$.  The new bound on
$E_{h+1}$ therefore produces
\[
 \begin{aligned}
 \binom n{h+1}|a_{h+1}-a|^{h+1}
 &\le |\bar f_{n-h-1}|+|E_{h+1}|+4^n\sum_{d=1}^h|E_d|\\
 &\le(1+H_n+h4^n)H_n^{2(h-1)}\epsilon
 \le H_n^{2h}\epsilon .
 \end{aligned}
\]
This establishes \eqref{eq:gw-completion-induction} and in particular proves 
\begin{equation}\label{eq:gw-completion-coefficients}
   |E_d|\le H_n^{2(q-1)}\epsilon \text{ for all } 1\le d\le q.
\end{equation} 

Since $z_{j_{n-1}}-a=0$,
every term of $E_d$ involving this entry is zero.  Thus $E_d$ is
also the $d$th elementary symmetric function of the $q+1$ entries
$z_k-a$ with $k\ne j_{n-1}$.  In particular,
\[
 E_{q+1}=\prod_{k\ne j_{n-1}}(z_k-a).
\]
For this choice of $a$ we have the polynomial identity
\[
 \prod_{k\ne j_{n-1}}(X-z_k)
 =\prod_{k\ne j_{n-1}}\bigl((X-a)-(z_k-a)\bigr)
 =\sum_{d=0}^{q+1}(-1)^d E_d(X-a)^{q+1-d}.
\]
If $a=0$, the factor indexed by $u$ gives $E_{q+1}=0$. Suppose $a\ne0$.
Then $j_1\ne j_{n-1}$, and both $j_1$ and $u$ contribute the factor $X$
on the left. Its constant and linear coefficients therefore vanish.
On the right, setting
the coefficient of $X$ equal to zero gives the identity below for
$(q+1)(-a)^q$. Further, setting the constant term equal to zero gives an expression for $E_{q+1}$.  Using $|a|\le1$ and the bound~\eqref{eq:gw-completion-coefficients} in these two identities gives
\[
 \begin{aligned}
 &(q+1)(-a)^q
 =-\sum_{d=1}^q(q+1-d)(-1)^dE_d(-a)^{q-d},\\
 &|a|^q
 \le \frac1{q+1}\sum_{d=1}^q(q+1-d)|E_d|
       \le \frac q2 H_n^{2(q-1)}\epsilon ,\\
 &|E_{q+1}|
 \le |a|^{q+1}+\sum_{d=1}^q |E_d|\,|a|^{q+1-d}
       \le \frac{3q}{2} H_n^{2(q-1)}\epsilon .
 \end{aligned}
\]
Combining the two cases with \eqref{eq:gw-completion-coefficients}
gives
\begin{equation}\label{eq:gw-completion-all-coefficients}
 |E_d|\le 2qH_n^{2(q-1)}\epsilon
 \qquad (1\le d\le q+1).
\end{equation}
Choose $k\ne j_{n-1}$ maximizing $|z_k-a|$.  As $-a$ is among the
remaining entries, $1=\max_l|z_l|\le2|z_k-a|$.  Setting $X=z_k$
in the displayed product identity and dividing by $|z_k-a|^{q+1}$ yields
\[
 1\le\sum_{d=1}^{q+1}\frac{|E_d|}{|z_k-a|^d}
 \le 2qH_n^{2(q-1)}\epsilon\sum_{d=1}^{q+1}2^d
 \le\delta_n^{-1}\epsilon .
\]
Here $\sum_{d=1}^{q+1}2^d=2^{q+2}-2\le4^{n-1}$, so
$\epsilon\ge\delta_n$.  For $0<\|\mathbf{y}\|_\infty\le1$, each form has
degree at most $q$, whereby multiplying \eqref{eq:gw-completion} by
$\|\mathbf{y}\|_\infty^q$ gives \eqref{eq:gw-completion-unweighted}.
\end{proof}

The estimate also restricts the limiting coefficient vectors in the
proof of Theorem~\ref{thm:analytic-equivalence}.  For $\nu-1\in\mathscr B_{\mathbf j}$ set
\[
 \mathscr F_\nu=\left\{\mathbf c\in\mathscr G_n:
     c_{n-h}=\frac{h+1}{\nu-h}\quad(1\le h<\nu)\right\}.
\]

\begin{corollary}\label{cor:first-change-prefix-centers}
\begin{enumerate}
\item It suffices to replace $\mathscr{G}_n$ by $\bigcup_{\nu-1\in\mathscr B_{\mathbf j}}\mathscr F_\nu$ in
Theorem~\ref{thm:analytic-equivalence}\textup{(ii)}, \textup{(iii)}.
Furthermore, the sets $\mathscr F_\nu$ are disjoint and
$|\mathscr F_\nu|=(q-\nu+1)!$.  If $\mathscr B_{\mathbf j}=\emptyset$, then Theorem~\ref{thm:analytic-equivalence}\textup{(i)} holds unconditionally.

\item Assume $C(n,w-1)$ and $|W_{\mathbf j}|=w\ge3$.  Then
$\nu=\min\{h\in\{2,\ldots,q\}:j_{n-h}\ne j_{n-1}\}$ exists, and
one may replace $\mathscr G^*_{n,\mathbf j}$ by $\mathscr F_\nu$ in
Theorem~\ref{thm:first-change-sliced-criterion}, taking neighborhoods
in $\mathcal H_{\mathbf j}$ as before.
\end{enumerate}
\end{corollary}
\begin{proof}
For either assertion, condition \textup{(i)} of
Theorem~\ref{thm:analytic-equivalence} implies the estimate near
the indicated centers by Theorem~\ref{thm:analytic-equivalence}\textup{(iv)}  and the
estimate implies uniform isolation.  It therefore suffices to show
that a nonzero common zero of $\bar f_2,\ldots,\bar f_{n-1}$
contradicts uniform isolation near the smaller set of centers.

\smallskip
\noindent\textup{(1)}
Let $p\ne\mathbf{0}$ be such a common zero, normalized as in the proof of
Theorem~\ref{thm:analytic-equivalence}. Let
\[
 a=a_1(p),\qquad E_0=1,\qquad
 E_d=e_d(z_k(p)-a:1\le k\le n).
\]
Suppose first that $a_h(p)=a$ for every $1\le h\le q$.  Then
$E_h=\bar f_{n-h}(p)=0$ for $1\le h\le q$.  For every
$h\in\mathscr B_{\mathbf j}$ we have $h+1\le q$, so that $\bar\gamma_{n-h}(p)=E_{h+1}-aE_h=0$. Thus all the forms in \eqref{eq:gw-completion} vanish at $p$,
contrary to Theorem~\ref{thm:gw-completion}.  Hence there exists
a least $s\in\{2,\ldots,q\}$ such that $a_s(p)\ne a$.
The index $j_{n-s}$ cannot have occurred in any smaller degree:
if it had, its value at $p$ would be $a$.  Thus
$s-1\in\mathscr B_{\mathbf j}$.  In particular, if
$\mathscr B_{\mathbf j}=\emptyset$, no such $p$ exists. For $1\le h<s$ we have $E_h=\bar f_{n-h}(p)=0$.  In degree $s$,
the translation identity \eqref{eq:elementary-translation} gives
\[
 0=\bar f_{n-s}(p)
   =\sum_{d=0}^{s}\binom{n-d}{s-d}
       (a-a_s(p))^{s-d}E_d
   =\binom ns(a-a_s(p))^s+E_s.
\]
Since $a_s(p)\ne a$, it follows that $E_s\ne0$. Now take the arc $(\mathbf{y}(t),\mathbf c(t))$ associated to $p$ in the
proof of Theorem~\ref{thm:analytic-equivalence},
\textup{(iii)}$\Rightarrow$\textup{(i)}.  For $1\le h<s$,
the form $\bar f_{n-h}$ has degree $h$ and selected coordinate
value $a$ at $p$.  Hence the least positive index $r_{n-h}$
in the limiting-coefficient formula of that proof is
\[
 r_{n-h}=\min\{r\ge1:E_{h+r}\ne0\}=s-h,
 \qquad
 c_{n-h}(0)=\frac{h+1}{r_{n-h}}=\frac{h+1}{s-h}.
\]
The limiting vector belongs to $\mathscr G_n$ by the construction,
and these equalities place it in $\mathscr F_s$.  Thus
\[
 \mathbf c(t)\to\mathbf c(0)\in\mathscr F_s,\quad
 \mathbf{y}(t)\to0,\quad
 \mathbf{y}(t)\ne0,\quad Q_{u,\mathbf c(t)}(\mathbf{y}(t))=0
 \quad(0<|t|\ll1).
\]
This contradicts uniform isolation near $\bigcup_{\nu-1\in\mathscr B_{\mathbf j}}\mathscr F_\nu$. Finally, $c_{n-1}=2/(s-1)$ on $\mathscr F_s$, so distinct values
of $s$ give disjoint sets.  On $\mathscr F_s$ the coordinates
$c_{n-1},\ldots,c_{n-s+1}$ are prescribed.  By
\eqref{eq:analytic-coefficient-grid}, each remaining coordinate
$c_i$, $2\le i\le n-s$, has $i-1$ independent choices.  Hence
\[
 |\mathscr F_s|=\prod_{i=2}^{n-s}(i-1)
               =(n-s-1)!=(q-s+1)!.
\]

\smallskip
\noindent\textup{(2)}
Assume the hypotheses of Theorem~\ref{thm:first-change-sliced-criterion},
and let $\nu$ be as in the statement.  The indices in degrees
$1,\ldots,\nu-1$ all equal $j_{n-1}$.  At degree $\nu$ the index
changes, and the proof of Theorem~\ref{thm:first-change-sliced-criterion}, by
Proposition~\ref{prop:regularity-with-difference}, gives $a_\nu(p)=z_{j_{n-\nu}}(p)\ne z_{j_{n-1}}(p)=a$. Thus $s=\nu$ in part~\textup{(1)}, and $\mathbf c(0)\in\mathscr F_\nu$.
In particular, taking $h=\nu-1$ in the coefficient formula gives
$c_{i_*}(0)=\nu$ for $i_*=n-\nu+1$.  The rescaling in the proof of
Theorem~\ref{thm:first-change-sliced-criterion} therefore produces an arc of nonzero zeros tending to the origin whose coefficients lie
in $\mathcal H_{\mathbf j}$ and tend to the same point
$\mathbf c(0)\in\mathscr F_\nu$.  This contradicts uniform
isolation in a neighborhood of $\mathscr F_\nu$ relative to
$\mathcal H_{\mathbf j}$, as required.
\end{proof}

\subsection{Uniform estimates at general coefficients}
\label{subsec:gw-generic}

We first compute the height of the ideal generated by any subset of
the universal equations.  The lower bound follows by applying
Theorem~\ref{thm:gw-completion} to a word chosen for that subset.

\begin{proposition}\label{prop:gw-universal-regularity}
Let $A=S_u[c_2,\ldots,c_{n-1}]=\mathbb C[y_1,\ldots,y_q,c_2,\ldots,c_{n-1}]$.
For every word $\mathbf j$, every $u\in U_{\mathbf j}$ and every
$T\subseteq\{1,\ldots,q\}$, one has
\begin{equation}\label{Eq:gw-universal-regularity}
 \operatorname{height}_A(Q_{n-h,u,\mathbf c}\mid h\in T)=|T|.
\end{equation}
In particular, every ordering of $Q_{2,u,\mathbf c},\ldots,
Q_{n-1,u,\mathbf c}$ is a regular sequence in $A$.
\end{proposition}
\begin{proof}
Fix a subset $S\subseteq\{1,\ldots,q\}$ and set $I_S=(\bar f_{n-h},\bar\gamma_{n-h}\mid h\in S)\subseteq S_u$. We first prove that $\dim S_u/I_S\le q-|S|$.  For this, choose a
word $\mathbf j'$ which agrees with $\mathbf j$ in the degrees
belonging to $S$, as follows.  Set $j'_1=j_1$, $j'_2=j_2$, and subsequently define
\[
 j'_i=
 \begin{cases}
   j_i,     & n-i\in S,\\
   j'_{i-1},& n-i\notin S
 \end{cases}
 \qquad(3\le i\le n-1).
\]
Every index occurring in $\mathbf j'$ occurs in $\mathbf j$, i.e., $W_{\mathbf{j}'}\subseteq W_{\mathbf{j}}$. Thus $u$ is still unused in $\mathbf{j}'$ i.e., $u\in U_{\mathbf{j}'}$ and $j'_1=j_1$ ensures that the analogues of $\bar f_i$ and $\bar \gamma_i$ for the word $\mathbf{j}'$, which we denote by  $\bar f'_i$, and $\bar\gamma'_i$ respectively, lie in the same ring $S_u$ as $\bar f_i$ and $\bar \gamma_i$. By \eqref{eq:coefficient-family}, $(\bar f'_{n-h},\bar\gamma'_{n-h})=(\bar f_{n-h},\bar\gamma_{n-h})$ for all $h\in S$. Moreover, if $1\le h<q$ and $h\notin S$, then
$j'_{n-h-1}=j'_{n-h}$ and thus, the selected index in degree $h+1$ is the same as in degree $h$.  Hence $\mathscr B_{\mathbf j'}\subseteq S$. It follows that the ideal
\[
\mathfrak{I}_S:= I_S+(\bar f'_{n-h}\mid 1\le h\le q,\ h\notin S)
\]
contains every $\bar f'_{n-h}$, $1\le h\le q$, as well as every
$\bar\gamma'_{n-h}$ with $h\in\mathscr B_{\mathbf j'}$.
If $\mathbf{y}\ne\mathbf{0}$ were a common zero of the ideal $\mathfrak{I}_S\subseteq S_u$, then
$\bar f'_{n-h}(\mathbf{y})=0$ for every $1\le h\le q$ and
$\bar\gamma'_{n-h}(\mathbf{y})=0$ for every $h\in\mathscr B_{\mathbf j'}$.
The left-hand side of \eqref{eq:gw-completion} for the word
$\mathbf j'$ and unused index $u$ would therefore be zero,
contrary to its lower bound $\delta_n>0$.  Since all the
generators are homogeneous of positive degree, the origin is
their only common zero. Thus the ideal $\mathfrak{I}_S/ I_S$ in $S_u/I_S$, generated by the image of 
$q-|S|$ additional forms, has radical $(y_1,\ldots,y_q)/I_S$ in $S_u/I_S$.  Since $I_S$ is homogeneous,
the dimension of $S_u/I_S$ equals its dimension at the origin.
Krull's height theorem (\cite[Lemmas~10.117.1 and~10.60.13]{stacks-project}) therefore yields
\[
 \dim S_u/I_S\le q-|S|
 \qquad\text{and hence}\qquad
 \operatorname{height} I_S\ge |S|
\]
.

Now fix $T\subseteq\{1,\ldots,q\}$ and let $\mathfrak P$ be a
prime ideal of $A=S_u[c_2,\ldots,c_{n-1}]$, containing $Q_{n-h,u,\mathbf c}$ for every
$h\in T$. Now set
\[
 \mathfrak p=\mathfrak P\cap S_u,\quad
 S=\{h\in T:\bar\gamma_{n-h}\in\mathfrak p\}.
\]
For $h\in S$, the identity
$Q_{n-h,u,\mathbf c}=\bar f_{n-h}+c_{n-h}\bar\gamma_{n-h}$
gives $\bar f_{n-h}\in\mathfrak p$.  Hence
$I_S\subseteq\mathfrak p$ and $\dim S_u/\mathfrak p\le q-|S|$.
For $h\in T\setminus S$, the same identity determines the image
of $c_{n-h}$ in $\operatorname{Frac}(A/\mathfrak P)$:
\[
 c_{n-h}=-\frac{\bar f_{n-h}}{\bar\gamma_{n-h}}
          \in\operatorname{Frac}(S_u/\mathfrak p),
\]
where the right-hand side denotes the quotient of the images
modulo $\mathfrak p$.  Thus at most $q-|T|+|S|$ of the coordinates of $\mathbf{c}=(c_2,\dots, c_{n-1})$ remain to generate $\operatorname{Frac}(A/\mathfrak P)$
over $\operatorname{Frac}(S_u/\mathfrak p)$.  Taking transcendence
degrees yields
\[
 \dim A/\mathfrak P
 \le \dim S_u/\mathfrak p+q-|T|+|S|
 \le (q-|S|)+(q-|T|+|S|)
 =2q-|T|.
\]
Since $\dim A=2q$, every such $\mathfrak P$ has height at least
$|T|$.  Conversely, the ideal generated by the $|T|$ indicated
forms is proper, since they all vanish at $\mathbf{y}=0$. Thus,  Krull's height theorem provides the reverse height inequality, thereby proving \eqref{Eq:gw-universal-regularity}.

Now consider some reordering $Q_{\sigma(2),u,\mathbf c},\ldots,Q_{\sigma(n-1),u,\mathbf c}$ of the sequence for some permutation $\sigma$. Let $\Sigma_m:=Q_{\sigma(2),u,\mathbf c},\ldots,Q_{\sigma(m),u,\mathbf c}$ be the initial partial sequence for any $2\leq m\leq n-1$. Applying \eqref{Eq:gw-universal-regularity} to the ideal generated by $\Sigma_m$, we have $\operatorname{height}_A(\Sigma_m)=m-1$ for all $2\leq m\leq n-1$. If the preceding $m-2$ terms form a regular sequence, every associated prime of their quotient has height $m-2$ in the polynomial ring $A$. Since $(\Sigma_m)$ has height $m-1$, the next term lies in none of these primes and is therefore a non-zero divisor on the quotient. Induction on $m$, starting with the empty sequence, establishes the regularity of the sequence $Q_{\sigma(2),u,\mathbf c},\ldots,Q_{\sigma(n-1),u,\mathbf c}$ in $A$.
\end{proof}
We can now prove the estimate \eqref{eq:finite-grid-estimate} for a general choice of $\mathbf{c}\in\mathbb{C}^q$.

\begin{theorem}\label{thm:gw-generic-estimate}
Let $n\ge3$, $q=n-2$, and $b_n=(n-1)!$.  There exists a nonzero
$\Delta_n\in\mathbb Q[c_2,\ldots,c_{n-1}]$ such that
$\Omega_n=\{\mathbf c\in\mathbb C^q:\Delta_n(\mathbf c)\ne0\}$
satisfies the following.
\begin{enumerate}[(i)]
\item For every word $\mathbf j$, every $u\in U_{\mathbf j}$ and every
$\mathbf c^0\in\Omega_n$,
\begin{equation}\label{Eq:fiberbound}
 1\le \dim_{\mathbb C}
 \frac{S_u}{(Q_{2,u,\mathbf c^0},\ldots,Q_{n-1,u,\mathbf c^0})}
 \le b_n.
\end{equation}
The displayed generators form a regular sequence in the local ring
$(S_u)_{(y_1,\ldots,y_q)}$.
\item For every compact $K\subset\Omega_n$ there exist $A_K,R_K>0$,
simultaneously for all words $\mathbf j$ and all $u\in U_{\mathbf j}$,
such that
\begin{equation}\label{eq:gw-generic-estimate}
 \|Q_{u,\mathbf c}(\mathbf{y})\|\ge A_K\|\mathbf{y}\|^{b_n}
 \qquad(\mathbf c\in K,\ \|\mathbf{y}\|<R_K).
\end{equation}
\end{enumerate}
Moreover, $\Omega_n\cap(\mathbb Q_{>0})^q\ne\varnothing$.
\end{theorem}
\begin{proof}
Fix a word $\mathbf j$ and an unused index $u\in U_{\mathbf j}$. Let $R=\mathbb Q[c_2,\ldots,c_{n-1}]$, $F=\operatorname{Frac}(R)$,
and
\[
 I_{\mathbf j,u}=(Q_{2,u,\mathbf c},\ldots,Q_{n-1,u,\mathbf c})
       \subset R[y_1,\ldots,y_q],\qquad
 B=F[y_1,\ldots,y_q]/I_{\mathbf j,u}F[y_1,\ldots,y_q].
\]
By Proposition~\ref{prop:gw-universal-regularity} and faithful
flatness of $\mathbb Q\subset\mathbb C$, the generators of
$I_{\mathbf j,u}$ form a regular sequence in $R[y_1,\ldots,y_q]$.
They remain regular after localization to $F[y_1,\ldots,y_q]$:
their ideal remains proper since each generator vanishes at $\mathbf{y}=0$.
Thus every prime containing this sequence of length $q$ has height
$q$ in $F[y_1,\ldots,y_q]$.  It follows that $B$ is a
zero-dimensional $F$-algebra of finite type, hence is
finite-dimensional over $F$. Choose an embedding $F=\mathbb Q(c_2,\ldots,c_{n-1})\hookrightarrow
\mathbb C$ by sending the $c_i$ to algebraically independent complex
numbers.  Since $B$ is finite-dimensional over $F$, the common
zeros of the resulting $q$ polynomials in $\mathbb C^q$ are isolated,
and the sum of their multiplicities is
$\dim_{\mathbb C}(B\otimes_F\mathbb C)=\dim_F B$.
The affine B\'ezout inequality
\cite[\S5.5, pp.~123--124]{Fulton1993Toric} therefore gives
\[
 1\le \dim_F B
 \le \prod_{i=2}^{n-1}\deg_{\mathbf{y}} Q_{i,u,\mathbf c}
 \le \prod_{i=2}^{n-1}(d_i+1)
 =(n-1)!=b_n.
\]
This establishes \eqref{Eq:fiberbound} on the geometric generic fiber. The same argument yields the same bound for any $\mathbf{c}^0\in \mathbb{C}^q$, for which $S_u/(Q_{2,u,\mathbf c^0},\ldots,Q_{n-1,u,\mathbf c^0})$ is zero-dimensional.

Let $N=\dim_F B$. For $1\le k\le q$, the $N+1$ residue classes
$1,y_k,\ldots,y_k^N$ are linearly dependent in $B$.  Thus there
is a monic polynomial $P_k(T)\in F[T]$ of degree at most $N$ with
$P_k(y_k)=0$ in $B$.  Evaluating at $\mathbf{y}=\mathbf{0}$ yields $P_k(0)=0$ and therefore we may write
\[
 P_k(T)=T^{s_k}H_k(T),\qquad H_k(0)\ne0,\qquad
 1\le s_k\le N\le b_n.
\]
By the definition of $B$, there are polynomials $A_{ki}\in F[\mathbf{y}]$
such that the following holds in $F[\mathbf{y}]$:
\begin{equation}\label{eq:gw-denominator-identities}
 y_k^{s_k}H_k(y_k)=P_k(y_k)
   =\sum_{i=2}^{n-1}A_{ki}(\mathbf{y})Q_{i,u,\mathbf c}(\mathbf{y})
 \qquad(1\le k\le q).
\end{equation}

Choose $\Delta_{\mathbf j,u}\in R\setminus\{0\}$ so that the
coefficients of all $P_k$ and $A_{ki}$, as well as the elements
$H_k(0)^{-1}$, belong to $R[1/\Delta_{\mathbf j,u}]$.
For every $\mathbf c^0\in \mathbb{C}^q$ satisfying
$\Delta_{\mathbf j,u}(\mathbf c^0)\ne0$,
\eqref{eq:gw-denominator-identities} specializes to a monic
relation for each $y_k$ in $S_u/(Q_{2,u,\mathbf c^0},\ldots,Q_{n-1,u,\mathbf c^0})$. This proves that this quotient is finite-dimensional over $\mathbb C$, since it is spanned by the residue classes of the monomials $y_1^{a_1}\cdots y_q^{a_q}$ with $0\le a_k<\deg P_k$. Thus, by the preceding B\'ezout bound, we obtain (the lower bound being a consequence of $Q_{i,u,\mathbf c}(\mathbf 0)=0$ for all $2\leq i\leq n-1$ for all $\mathbf{c}\in \mathbb{C}^q$):
\[1\le
 \dim_{\mathbb C}
 \frac{S_u}{(Q_{2,u,\mathbf c^0},\ldots,Q_{n-1,u,\mathbf c^0})}
 \le b_n\]
for all $\mathbf c^0\in \mathbb{C}^q$ satisfying
$\Delta_{\mathbf j,u}(\mathbf c^0)\ne0$. Consequently, $Q_{2,u,\mathbf c^0},\ldots,Q_{n-1,u,\mathbf c^0}$ is a system of parameters, and hence a regular sequence in the
$q$-dimensional regular local ring $(S_u)_{(y_1,\ldots,y_q)}$ by \cite[Lemma~10.104.2]{stacks-project}. This establishes Theorem~\ref{thm:gw-generic-estimate}\textup{(i)}.

It remains to prove the uniform estimate in \textup{(ii)}.  Fix a nonempty compact
$K\subset\{\Delta_{\mathbf j,u}\ne0\}$, and write
$H_k(T,\mathbf c)$ and $A_{ki}(\mathbf{y},\mathbf c)$ for the specializations of the polynomials in \eqref{eq:gw-denominator-identities} at some $\mathbf{c}\in\mathbb{C}^q$.  Since $H_k(0,\mathbf c)\ne0$ on $K$, compactness
and continuity give $a,M>0$ and $0<R_K\le1$ such that
\[
 |H_k(y_k,\mathbf c)|\ge a,\qquad
 |A_{ki}(\mathbf{y},\mathbf c)|\le M
 \qquad(\mathbf c\in K,\ \|\mathbf{y}\|<R_K)
\]
for all $k,i$.  Choose $k$ such that $|y_k|=\|\mathbf{y}\|_\infty$.  Since
$|y_k|\le\|\mathbf{y}\|<R_K\le1$ and $s_k\le b_n$, we have
$|y_k|^{b_n}\le |y_k|^{s_k}$.  Using
$\|\mathbf{y}\|\le\sqrt q\,\|\mathbf{y}\|_\infty$ and
\eqref{eq:gw-denominator-identities}, we obtain
\[
 \begin{aligned}
 q^{-b_n/2}\|\mathbf{y}\|^{b_n}\le |y_k|^{b_n}
 \le |y_k|^{s_k} =\frac{\bigl|\sum_{i=2}^{n-1}
       A_{ki}(\mathbf{y},\mathbf c)Q_{i,u,\mathbf c}(\mathbf{y})\bigr|}
        {|H_k(y_k,\mathbf c)|} \le \frac{M}{a}\sum_{i=2}^{n-1}|Q_{i,u,\mathbf c}(\mathbf{y})|
 \le \frac{\sqrt q\,M}{a}\,\|Q_{u,\mathbf c}(\mathbf{y})\|,
 \end{aligned}
\]
where the last inequality is Cauchy--Schwarz, since there are
$n-2=q$ summands.
This proves the required estimate for the fixed pair $(\mathbf j,u)$. Finally, let $\Delta_n=\prod_{(\mathbf j,u)}\Delta_{\mathbf j,u}$, where the product ranges over the finitely many words of length $n-1$ and their unused indices.  For compact $K\subset\Omega_n$,
taking the minima of the lower-bound constants and of the radii
over these pairs gives \eqref{eq:gw-generic-estimate}.  Since
$(\mathbb Q_{>0})^q$ is Zariski dense in $\mathbb A^q_{\mathbb Q}$
and $\Delta_n\ne0$, there exists
$\mathbf c^0\in(\mathbb Q_{>0})^q$ with $\Delta_n(\mathbf c^0)\ne0$.
\end{proof}

\subsection{A common coefficient open set with maximal grid intersection}
\label{subsec:grid-count}

Fix $n\ge3$, and let $q=n-2$, $b_n=(n-1)!$ and
$R=\mathbb Q[c_2,\ldots,c_{n-1}]$ as before.
\begin{definition}\label{Def:admissible}
A Zariski open set $O\subseteq\mathbb C^q$ is \emph{admissible} if both of the following hold:
\begin{enumerate}[(i)]
\item For every $\mathbf c^0\in O$, every word $\mathbf j$ and every
$u\in U_{\mathbf j}$, the quotient
$S_u/(Q_{2,u,\mathbf c^0},\ldots,Q_{n-1,u,\mathbf c^0})$
is finite-dimensional over $\mathbb C$.
\item For every compact $K\subset O$, there exist $A_K,R_K>0$ such
that \eqref{eq:gw-generic-estimate} holds for all $\mathbf c\in K$,
all $\|\mathbf{y}\|<R_K$, and all pairs $(\mathbf j,u)$.
\end{enumerate}
\end{definition}
Theorem~\ref{thm:gw-generic-estimate} gives a nonempty admissible
principal open set.  We now determine the largest admissible open,
and a principal open with the same real points.
For each $(\mathbf j,u)$ retain $I_{\mathbf j,u}\subset R[\mathbf{y}]$ from the
proof of that theorem, put $\mathfrak m=(y_1,\ldots,y_q)$, and
write $\operatorname{ev}_{\mathbf0}:R[\mathbf{y}]\to R$ for evaluation at $\mathbf{y}=0$.
We take zero sets over $\mathbb C$.  For ideals $J,L$ we use the saturation
$J:L^\infty=\bigcup_{k\ge1}(J:L^k)$, so that
$V(J:L^\infty)=\overline{V(J)\setminus V(L)}^{\,\mathrm{Zar}}$.
Indeed, on the principal open sets defined by generators of $L$ this
is localization followed by contraction, and one takes the union of
the resulting closed sets.  These ideal quotients, finite intersections
and contractions commute with extension from $\mathbb Q$ to $\mathbb C$
by flatness.

\begin{lemma}\label{lem:grid-optimal-local}
Fix $(\mathbf j,u)$, and let $\mathfrak b_{\mathbf j,u}^{\,0}
 =\operatorname{ev}_\mathbf{0}(I_{\mathbf j,u}:\mathfrak m^{b_n})\subset R$.
For $\mathbf c^0\in\mathbb C^q$, the following are equivalent:
\begin{enumerate}[(i)]
\item $\mathbf c^0\notin V(\mathfrak b_{\mathbf j,u}^{\,0})$.
\item There exist $A,\varepsilon>0$ such that
\[
 \|Q_{u,\mathbf c}(\mathbf{y})\|\ge A\|\mathbf{y}\|^{b_n}
 \qquad
 (\|\mathbf{y}\|<\varepsilon,\ \|\mathbf c-\mathbf c^0\|<\varepsilon).
\]
\item The estimate in \textup{(ii)} holds with some exponent $B>0$
in place of $b_n$.
\end{enumerate}
\end{lemma}

\begin{proof}
By Proposition~\ref{prop:gw-universal-regularity} and faithful
flatness, $I_{\mathbf j,u}$ is generated by the regular sequence $Q_{2,u,\mathbf c},\ldots,Q_{n-1,u,\mathbf c}$ of
length $q$ in $R[\mathbf{y}]$.  Hence all its associated primes have height
$q$.  Since $I_{\mathbf j,u}\subseteq\mathfrak m$ and
$\operatorname{height}\mathfrak m=q$, the only associated prime
containing $\mathfrak m$ is $\mathfrak m$ itself.
In a minimal primary decomposition write
$I_{\mathbf j,u}=J_0\cap J'$, where $J_0$ is $\mathfrak m$-primary
and $J'$ is the intersection of the other components (or $R[\mathbf{y}]$
if there are none). The ring $(R[\mathbf{y}]/I_{\mathbf j,u})_{\mathfrak m}$ is the local ring at $\mathbf{y}=\mathbf{0}$ of the generic fiber over $F=\operatorname{Frac}(R)$ and is thus an Artinian $F$-algebra. Its length is bounded above by the $F$-dimension of the generic fiber algebra,
which is at most $b_n$ by the proof of
Theorem~\ref{thm:gw-generic-estimate}.  In an Artinian local ring of
length $\ell$, the $\ell^{th}$ power of the maximal ideal is zero, whereby
\[
 \mathfrak m^{b_n}\subseteq
 I_{\mathbf j,u}R[\mathbf{y}]_{\mathfrak m}\cap R[\mathbf{y}]=J_0.
\]
Let $J$ be a primary component occurring in $J'$.  The primes
$\sqrt J$ and $\mathfrak m$ are distinct and both have height $q$.
Hence $\mathfrak m\not\subseteq\sqrt J$, so $y_k\notin\sqrt J$
for some $k$.  If $f\in J:\mathfrak m^s$, then $fy_k^s\in J$.
Since $y_k^s\notin\sqrt J$ and $J$ is primary, this implies
$f\in J$.  Thus $J:\mathfrak m^s=J$ for every $s\ge1$.
Intersecting over these components gives
$J':\mathfrak m^s=J'$. On the other hand, for $s\ge b_n$ we have
$\mathfrak m^s\subseteq\mathfrak m^{b_n}\subseteq J_0$, whence
$J_0:\mathfrak m^s=R[\mathbf{y}]$.  Therefore
\[
 I_{\mathbf j,u}:\mathfrak m^s
 =(J_0:\mathfrak m^s)\cap(J':\mathfrak m^s)
 =R[\mathbf{y}]\cap J'=J'
 \qquad(s\ge b_n).
\]
Thus, by the definition of saturation,
\begin{equation}\label{eq:grid-bounded-saturation}
 I_{\mathbf j,u}:\mathfrak m^{b_n}
 =J'=I_{\mathbf j,u}:\mathfrak m^\infty .
\end{equation}

If \textup{(i)} holds, choose $g\in I_{\mathbf j,u}:\mathfrak m^{b_n}$
with $g(0,\mathbf c^0)\ne0$.  For each $k$ there are $A_{ki}\in R[\mathbf{y}]$
such that
\[
 g\,y_k^{b_n}=\sum_{i=2}^{n-1}A_{ki}Q_{i,u,\mathbf c}.
\]
On a sufficiently small product of balls about $(0,\mathbf c^0)$
all $|A_{ki}/g|$ are bounded by a common $M>0$.  Choosing
$|y_k|=\|\mathbf{y}\|_\infty$ gives
\[
 q^{-b_n/2}\|\mathbf{y}\|^{b_n}
 \le |y_k|^{b_n}
 \le M\sum_{i=2}^{n-1}|Q_{i,u,\mathbf c}(\mathbf{y})|
 \le \sqrt q\,M\|Q_{u,\mathbf c}(\mathbf{y})\|.
\]
This proves \textup{(ii)}, which implies \textup{(iii)}. Conversely, if $\mathbf c^0\in V(\mathfrak b_{\mathbf j,u}^{\,0})$,
then $(0,\mathbf c^0)\in V(J')$ by \eqref{eq:grid-bounded-saturation}
and the definition of evaluation.  By saturation, $V(J')$ is the
Zariski closure of $V(I_{\mathbf j,u})\setminus V(\mathfrak m)$.
For a constructible subset of complex affine space its Zariski
and Euclidean closures agree by \cite[Aside~9.8]{milneAG}.
Hence there are $(\mathbf{y}^{(s)},\mathbf c^{(s)})\to(0,\mathbf c^0)$ with
$\mathbf{y}^{(s)}\ne0$ and $Q_{u,\mathbf c^{(s)}}(\mathbf{y}^{(s)})=0$.
These zeros contradict \textup{(iii)} for every $B>0$.
\end{proof}

\begin{lemma}\label{lem:grid-optimal-fibers}
For each $(\mathbf j,u)$ there is an explicitly determined ideal
$\mathfrak b_{\mathbf j,u}^{\,\mathrm{fib}}\subset R$ such that
\[
 V(\mathfrak b_{\mathbf j,u}^{\,\mathrm{fib}})
 =\overline{\{\mathbf c^0\in\mathbb C^q:
                \dim V(Q_{u,\mathbf c^0})>0\}}^{\,\mathrm{Zar}} .
\]
\end{lemma}
\begin{proof}
Every finite fiber has length at most $b_n$, by the affine
B\'ezout bound in the proof of Theorem~\ref{thm:gw-generic-estimate},
and contains the origin.  A positive-dimensional complex affine variety
has infinitely many points.  Thus a fiber has positive dimension
if and only if it contains $b_n$ distinct nonzero points.

Let $\mathbf{y}^{(0)}=0$ and introduce $b_n$ independent $q$-tuples of
variables $\mathbf{y}^{(1)},\ldots,\mathbf{y}^{(b_n)}$.  In
$R[\mathbf{y}^{(1)},\ldots,\mathbf{y}^{(b_n)}]$ put
\[
 \begin{aligned}
 J&=\bigl(Q_{i,u,\mathbf c}(\mathbf{y}^{(s)})
       \mid 2\le i\le n-1,\ 1\le s\le b_n\bigr),\\
 D&=\prod_{0\le s<t\le b_n}
    (y^{(s)}_1-y^{(t)}_1,\ldots,y^{(s)}_q-y^{(t)}_q).
 \end{aligned}
\]
The zero set of $D$ is the union of the diagonals
$\mathbf{y}^{(s)}=\mathbf{y}^{(t)}$, including $\mathbf{y}^{(t)}=0$ for $s=0$.  The image of
$V(J)\setminus V(D)$ in coefficient space therefore consists precisely
of the parameters with positive-dimensional fiber.  Taking its closure
amounts to saturating to remove the diagonals and then eliminating the
$\mathbf{y}^{(s)}$ \cite[Lemma~29.6.3]{stacks-project}.  Thus the required ideal is
\[
 \mathfrak b_{\mathbf j,u}^{\,\mathrm{fib}}=(J:D^\infty)\cap R .
\]
\end{proof}

The following result identifies the maximal admissible open subset of $\mathbb{C}^q$. Furthermore, we also obtain explicit choices for the polynomial $\Delta_n\in \mathbb{Q}[c_2,\dots, c_{n-1}]$ in Theorem~\ref{thm:gw-generic-estimate} such that the corresponding $\Omega_n\subseteq \mathbb{C}^q$ contain the maximal number of points of $\Gamma= \mathscr G_n$, $\mathscr G^*_{n,\mathbf j}$ and, $\mathscr F_\nu$ among all admissible open subsets of $\mathbb{C}^q$.
\begin{theorem}\label{thm:grid-optimal-open}
Let
\[
 \mathfrak b_n=\bigcap_{(\mathbf j,u)}
   \bigl(\mathfrak b_{\mathbf j,u}^{\,0}
        \cap\mathfrak b_{\mathbf j,u}^{\,\mathrm{fib}}\bigr)\subset R,
\]
where the intersection runs over all words $\mathbf{j}=(j_1,\dots, j_{n-1})\in[n]^{n-1}$ and all unused indices $u\in U_\mathbf{j}$.  Then $U=\mathbb C^q\setminus V(\mathfrak b_n)$
is the largest admissible Zariski open set.  For any generators
$g_1,\ldots,g_s\in R$ of $\mathfrak b_n$, let
\begin{equation}\label{eq:grid-optimal-denominator}
 \Delta_n=\sum_{a=1}^s g_a^2,\qquad
 \Omega_n=\{\mathbf c\in\mathbb C^q:\Delta_n(\mathbf c)\ne0\}.
\end{equation}
Then $\Omega_n$ is admissible and $\Omega_n\cap\mathbb R^q=U\cap\mathbb R^q$.
In particular, for every admissible $O$, one has
$O\cap\Gamma\subseteq\Omega_n\cap\Gamma$ for $\Gamma=\mathscr G_n$,
and also for each $\mathscr G^*_{n,\mathbf j}$ and $\mathscr F_\nu$
defined in Theorem~\ref{thm:first-change-sliced-criterion} and
Corollary~\ref{cor:first-change-prefix-centers}.
\end{theorem}
\begin{proof}
Let $U=\mathbb C^q\setminus V(\mathfrak b_n)$ and $O\subseteq \mathbb{C}^q$ be an arbitrary admissible open subset. For $\mathbf c^0\in O$, take a closed Euclidean ball of positive radius cenetered at $\mathbf c^0$, which is contained in $O$.  Applying Definition~\ref{Def:admissible}\textup{(ii)} to this compact ball and then Lemma~\ref{lem:grid-optimal-local} implies
$\mathbf c^0\notin V(\mathfrak b_{\mathbf j,u}^{\,0})$ for every
pair $(\mathbf{j}, u)$.  By Definition~\ref{Def:admissible}\textup{(i)} and Lemma~\ref{lem:grid-optimal-fibers}, $O\cap V(\mathfrak b_{\mathbf j,u}^{\,\mathrm{fib}})=\emptyset$.  Since
\[
 V(\mathfrak b_n)=\bigcup_{(\mathbf j,u)}
  \bigl(V(\mathfrak b_{\mathbf j,u}^{\,0})\cup
        V(\mathfrak b_{\mathbf j,u}^{\,\mathrm{fib}})\bigr),
\]
we obtain $O\subseteq U$. Conversely, $\dim V(Q_{2,u,\mathbf c^0},\ldots,Q_{n-1,u,\mathbf c^0})=0$ for all $\mathbf{c}^0\in U$ by Lemma~\ref{lem:grid-optimal-fibers}. Now let $K\subset U$ be an arbitrary compact subset and fix a pair $(\mathbf j,u)$.  By Lemma~\ref{lem:grid-optimal-local}, each
$\mathbf c^0\in K$ has an open neighborhood in $U$ on which the
estimate of Lemma~\ref{lem:grid-optimal-local}\textup{(ii)} holds with exponent $b_n$. These open neighborhoods form an open cover of $K$, which admits a finite subcover due to compactness of $K$. Taking the minima of the associated constants $A, \epsilon$ over the finitely many open neighborhoods covering $K$ and then over the finitely many pairs $(\mathbf j,u)$ proves the estimate \eqref{eq:gw-generic-estimate} on $K$. Since $K\subseteq U$ is an arbitrary compact subset, this proves that $U$ is admissible.

Finally, by \eqref{eq:grid-optimal-denominator}, $\Delta_n\in\mathfrak b_n$, so $\Omega_n\subseteq U$.
For $\mathbf c\in\mathbb R^q$,
\[
 \Delta_n(\mathbf c)=0
 \quad\Longleftrightarrow\quad
 g_a(\mathbf c)=0\ \text{for every }a
 \quad\Longleftrightarrow\quad \mathbf c\notin U.
\]
Thus $\Omega_n\cap\mathbb R^q=U\cap\mathbb R^q$.  By
Theorem~\ref{thm:gw-generic-estimate}, there exists an admissible open set intersecting $(\mathbb{Q}_{>0})^q$. By the above argument, this point belongs to $U$ and hence to $\Omega_n$ as well, whereby $\Delta_n\ne0$. The final grid inclusions $O\cap\Gamma\subseteq\Omega_n\cap\Gamma$ for $\Gamma=\mathscr G^*_{n,\mathbf j},\  \mathscr G_n, \ \mathscr F_\nu$ follow since they are subsets of $\mathbb{Q}^q$.
\end{proof}

\begin{corollary}\label{prop:grid-open-count}
Let $n\ge4$, and fix a word $\mathbf j=(j_1,\ldots,j_{n-1})$
such that $|\{j_2,\ldots,j_{n-1}\}|\geq 2$. Define $\nu=\min\{h\in\{2,\ldots,q\}:j_{n-h}\ne j_{n-1}\}$ and $i_*=n-\nu+1$. Furthermore, let $\Gamma=\mathscr G^*_{n,\mathbf j}=\mathscr G_n\cap\{c_{i_*}=\nu\}$, and let $\Omega_n$ be the
principal open set defined in \eqref{eq:grid-optimal-denominator}. Then $|\Gamma|=q!/(n-\nu)$ and the vanishing ideal of $\Gamma$ in $R$ is
\[
 J_*=\left(c_{i_*}-\nu,\ \prod_{r=1}^{i-1}\left(c_i-\frac{d_i+1}{r}\right) \text{ for all } 2\le i\le n-1,\ i\ne i_*
 \right)\subseteq R,
\]
 Moreover,
\begin{equation}\label{eq:exact-grid-open-count}
 \max_{{\rm admissible}\ O\subseteq\mathbb{C}^q}\frac{|O\cap\Gamma|}{|\Gamma|}
 =\frac{|\Omega_n\cap\Gamma|}{|\Gamma|}
 =1-\frac{\dim_{\mathbb Q}\bigl(R/(J_*, \ \mathfrak b_n)\bigr)}
             {|\Gamma|}.
\end{equation}
\end{corollary}

\begin{proof}
By \eqref{eq:analytic-coefficient-grid}, the $i^{th}$ coordinate
of points in $\mathscr G_n$ admits the $i-1$ distinct rational values
$(d_i+1)/r$, $1\le r\le i-1$.  Since $d_{i_*}+1=\nu$, the
condition $c_{i_*}=\nu$ uniquely determines the $i_*^{th}$ coordinate for points in $\Gamma=\mathscr G^*_{n,\mathbf j}=\mathscr G_n\cap\{c_{i_*}=\nu\}$. All other coordinates vary independently. Thus, the explicit description of the points of $\Gamma$ yields its size and the description of the vanishing ideal $J_\ast\subseteq R$ of $\Gamma $. Since $\Gamma\subset\mathbb Q^q$, Theorem~\ref{thm:grid-optimal-open} yields $\Gamma\setminus\Omega_n=\Gamma\cap V(\mathfrak b_n)$ and $O\cap\Gamma\subseteq\Omega_n\cap\Gamma$ for every admissible open $O\subseteq\mathbb{C}^q$. For any ideal $I\subseteq R$, it follows that $\dim_{\mathbb Q}\bigl(R/(J_*+I)\bigr)=|\Gamma\cap V(I)|$. In particular, taking $I=\mathfrak b_n$ proves \eqref{eq:exact-grid-open-count}, since $\Omega_n$ itself is admissible.
\end{proof}

For $\Gamma=\mathscr F_\nu$ the same calculation applies after
adjoining its defining linear equations to $J_*$ and replacing
$|\Gamma|$ by $(q-\nu+1)!$.  

\printbibliography
\end{document}